\documentclass[a4paper, 11pt]{amsart}
\usepackage{graphicx} % Required for inserting images
\usepackage{amsmath,color,soulutf8,longtable,colortbl,setspace,ifthen,xspace,url,pdflscape,tikz,pgfplots,amsthm}
\usepackage[utf8]{inputenc}
\usepackage{lmodern} 
\usepackage[T1]{fontenc}
\usepackage{amssymb}
\usepackage{amsfonts}
\usepackage{tikz,tikz-cd}
\usepackage{adjustbox}
\usepackage{xcolor}
\usepackage{hyperref}
\usepackage{caption}  % Hyperlien vers la figure plutôt que son titre.
\usepackage{pifont}
\usepackage{esint}
\usepackage{adjustbox}
\usepackage{geometry}
\usepackage{enumerate}
\usepackage{enumitem}
\usepackage{mathrsfs} 
\usepackage{xcolor}
\usepackage{ stmaryrd }
\tikzset{%
    symbol/.style={%
        draw=none,
        every to/.append style={%
            edge node={node [sloped, allow upside down, auto=false]{$#1$}}}
    }
}

\theoremstyle{definition}
\newtheorem{definition}{Definition}
\theoremstyle{definition}
\newtheorem{theorem}[definition]{Theorem}
\theoremstyle{definition}

\newtheorem{lemma}[definition]{Lemma}
\theoremstyle{definition}
\newtheorem{example}[definition]{Example}
 
\theoremstyle{definition}
\newtheorem{remarque}[definition]{Remark}
\theoremstyle{construction}

\theoremstyle{definition}
\newtheorem*{demo}{Proof :}
\theoremstyle{definition}
\renewenvironment{demo}{{\bfseries Proof :}}{\hfill\qedsymbol}
\newtheorem{proposition}[definition]{Proposition}
\theoremstyle{esquisse}

\theoremstyle{definition}

\theoremstyle{definition}

\newcommand{\A}{\mathcal{A}}
\newcommand{\p}{\mathcal{P}}
\newcommand{\Hom}{\underline{Hom}_\A}
\newcommand{\Homgr}{\underline{Hom}_{gr\A}}
\newcommand{\Homch}{\underline{Hom}_{Ch(\A)}}
\mathchardef\mhyphen="2D

\title{Higher crossed modules of algebras over an operad}
\author{Clovis Chabertier} 
\begin{document}
\maketitle
\tableofcontents
\section*{Introduction}

Crossed modules have been studied in various contexts for a long time in algebraic topology, beginning with the work of Whitehead in \cite{bams/1183513797} to understand pointed homotopy $2\mhyphen$types. Given a pair of pointed spaces $A \subset X$, the crossed module associated to it is the following. The group $\pi_1(A,\ast)$ acts by automorphisms on $\pi_2(X,A,\ast)$ and the connecting group morphism $\partial : \pi_2(X,A,\ast)\rightarrow \pi_1(A,\ast)$ satisfies equivariance with respect to this action, and also satisfies the Peiffer condition : ${}^{\partial x}y=xyx^{-1}$.
Loday in \cite{loday1982spaces} shows that crossed modules of groups can be understood in many ways : as groups internal to the category of small categories, as simplicial groups with Moore complex of length $1$, as $\mathrm{Cat}^1\mhyphen$groups, or as categories internal to the category of groups. These last two descriptions allow him to generalize crossed modules of groups to higher versions, namely $\mathrm{Cat}^n\mhyphen$groups. Morally such an object is an $n\mhyphen$fold category internal to groups. Thus higher versions of crossed modules, say $n\mhyphen$crossed modules, can be inductively defined as crossed modules internal to the category of $(n-1)\mhyphen$crossed modules. This procedure of internalization yields an equivalence of categories between the category of $n\mhyphen$crossed modules of groups and $\mathrm{Cat}^n\mhyphen$groups. The other idea of Loday is that one can associate to an $n\mhyphen$crossed module $\mathcal{L}=\mathbb{H}\overset{d}{\rightarrow} \mathbb{G}$ a space $B\mathcal{L}$ such that the canonical sequence $B\mathbb{H}\rightarrow B\mathbb{G}\rightarrow B\mathcal{L}$ is a homotopy fiber sequence of pointed spaces. Using these ideas, he was able to prove that $n\mhyphen$crossed modules are models for pointed homotopy $(n+1)\mhyphen$types. Later on, several authors studied crossed modules in algebraic contexts. For example Ellis in \cite{ELLIS1988277} studies $n\mhyphen$crossed cubes of associative/commutative/Lie/etc.. algebras and proves that such objects are equivalently defined as $n\mhyphen$fold categories internal to algebras of the associated type. In 2003, Janelidze in \cite{janelidze2003internal} gives a general framework in which crossed modules can be defined. This framework is the one of semi-abelian categories, and he proves that the category of crossed modules internal to a semi-abelian category $\mathcal{C}$ is equivalent to the category $\mathrm{Cat}(\mathcal{C})$ of internal categories to $\mathcal{C}$. As the category of categories internal to a semi-abelian category is again semi-abelian, we can iterate this construction and prove that $n\mhyphen$fold crossed modules internal to $\mathcal{C}$ are equivalent to $n\mhyphen$fold categories internal to $\mathcal{C}$. Many algebraic categories are known to be semi-abelian : the categories of groups, non-unital rings, associative or commutative algebras, Lie algebras, etc ... And it is folklore that the category of algebras over a symmetric reduced algebraic (see \cite{loday2012algebraic}) operad $\p$ is semi-abelian. 

However, the approach of Janelidze and Ellis becomes intricate when one wants to study higher crossed modules of $\p\mhyphen$algebras. Indeed, the approach of Ellis requires quite a few axioms, and the one of Janelidze requires to compute coproducts of $\p\mhyphen$algebras, which is hardly computable in pratice. Moreover, (higher) crossed modules of $\p\mhyphen$algebras (or of any kind) must have an homotopical flavour linked to the homotopy theory of $\p\mhyphen$algebras, namely differential graded $\p\mhyphen$algebras.

%However the definitions of Janelidze and Ellis are not so well adapted to the case of algebras over an operad. First of all, their definitions are a bit involved, Indeed they both required lots of axioms. For example, a crossed module of associative algebras is the data of a morphism of algebras $d:A\rightarrow B$, an internal action of $B$ on $A$ such that $d$ is equivariant and satisfies a Peiffer condition. When one wants to go to higher crossed modules, this becomes complicated to have a concrete grasp on it. Second, the homotopical properties, and especially the link between (higher) crossed modules of algebras and the homotopy theory of algebras over an operad is not entirely clear in their approach. 

In \cite{bauesminianrichter}, Baues, Minian and Richter propose another approach for crossed modules of $\p\mhyphen$algebras to obtain an operadic version of Hochschild cohomology. A slight variation of their definition is the starting point of a current work of Leray, Rivière and Wagemann in \cite{lrw}. Namely, a crossed module of $\mathcal{P}\mhyphen$algebras is a $\p\mhyphen$algebra structure on a chain complex $...0\rightarrow V_1 \overset{d}{\rightarrow} V_0$ concentrated in degrees $0$ and $1$. This last approach to crossed modules of $\p\mhyphen$algebras is the one we choose to work with in the present article.

We also would like to mention the work of Lack and Paoli in \cite{lackpaoli}, where they study the interplay between finite limit structures and operadic structures. In particular they prove that for any operad $\p$ in a closed symmetric monoidal category $\mathcal{V}$, there exists an operad $\p'$ in the category $\mathrm{Cat}(\mathcal{V})$, such that there is an equivalence of categories $\mathrm{Cat}(\p\mhyphen alg)\simeq \p'\mhyphen alg$. In fact they prove this statement for other finite limit theories than $Cat$ under some mild assumptions. This viewpoint, as soon as crossed modules of $\p\mhyphen$algebras and categories internal to $\p\mhyphen$algebras are identified, allows an operadic treatment, and all the machinery operads comes with, of crossed modules. This machinery includes Koszul duality, model structure for operads, and more generally the tools of homotopical algebra. The approach of Leray-Rivière-Wagemann is kind of orthogonal to this one as, instead of replacing the operad $\p$ by another one $\p'$, they replace the category where the operad $\p$ acts. 

The main goal of this paper is to study the approach, and its apparent easiness, of (higher) crossed modules given by Leray-Rivière-Wagemann in \cite{lrw} and highlight its link with differential graded $\p\mhyphen$algebras. We prove that this new approach is equivalent to previous ones, say of Janelidze or Ellis, for the category of $\p\mhyphen$algebras. To do so, we exhibit a codescent (or derived operations in reference to Kosmann-Schwarzbach \cite{Kosmann_Schwarzbach_2004}) phenomena that happens at the level of $\p\mhyphen$algebra structures. Indeed, given a (global) $\p\mhyphen$algebra structure on a chain complex $V=V_1 \rightarrow V_0$, we define $\p\mhyphen$algebra structures on $V_1$ and $V_0$ such that $V_1\rightarrow V_0$ is a crossed module of $\p\mhyphen$algebras in the sense of Janelidze. In order to define such a $\p\mhyphen$algebra structure on $V_1$ we exhibit an operad morphism $\rtimes : End^\A(V) \rightarrow End(V_1\oplus V_0)$ (see Proposition \ref{codes}). This operad morphism satisfies the following property : given a $\p\mhyphen$algebra structure on $V$, there exists a $\p\mhyphen$algebra structure on $V_1\oplus V_0$ such that the second projection $p_2:V_1\oplus V_0 \rightarrow V_0$ is a $\p\mhyphen$algebra morphism. It implies that $V_1$ is a $\p\mhyphen$algebra, and as a $\p\mhyphen$algebra, $V_1\oplus V_0$ is a semi-direct product $V_1 \rtimes V_0$. Moreover, this construction of the semi-direct product is somehow universal, in the sense that it is constructed at the operadic level, independantly of $\p$, and not in the category of $\p\mhyphen$algebra itself, as it used to be. 
\subsection*{Organisation of the paper}
The first part of this paper is devoted to a reformulation of the work of Janelidze in the context where the semi-abelian category is the category of $\p\mhyphen$algebras. This reformulation behaves nicely in the operadic context, in particular we extract, in Definition \ref{Peiffer}, some Peiffer relations, which are well known in the Lie and associative algebra cases \cite{baez2003higher, baues2002crossed,wagemann2021crossed}. We also provide an explicit description of $\mathrm{Cat}^1\mhyphen\p\mhyphen$algebras. This notion finds its roots in the following observation, due to Grothendieck \cite{grothendieck1968categories}. A reflexive graph in an abelian category $\A$ admits a unique composition which promotes it to an internal category to $\A$. As a consequence, given a category $\mathcal{C}$ endowed with a faithful right adjoint $U:\mathcal{C}\rightarrow \A$, then on reflexive graphs internal to $\mathcal{C}$, there exists at most one composition which turns it into an internal category to $\mathcal{C}$. Therefore the existence of such a composition becomes a property of this graph and no more a data on the graph itself. If $\mathcal{C}$ is the category of $\p\mhyphen$algebras in $\A$, then the faithful right adjoint is given by the forgetful functor $U:\mathcal{C} \rightarrow \A$. The $\mathrm{Cat}^1\mhyphen$condition of Definition \ref{cat^1} is an explicit necessary and sufficient condition for a reflexive graph in $\mathcal{C}$ to admit an internal category structure. The operadic treatment of crossed modules allows us to find back the explicit notion of crossed modules of Lie algebras \cite{baez2003higher}, associative algebras \cite{baues2002crossed}, Leibniz algebras and diassociative algebras \cite{casas2015crossedmoduleslieleibniz}.

The second part of this paper is devoted to the study of crossed modules in the sense of Leray-Rivière-Wagemann, and its link with other notions of crossed modules of $\p\mhyphen$algebras. We prove that both approaches to crossed modules of $\p\mhyphen$algebras are equivalent. This new approach is thought as a global approach, in the sense that a crossed module is a \textit{global} $\p\mhyphen$algebra structure on a chain complex. This approach is therefore linked to the homotopy theory of $\p\mhyphen$algebras. Moreover, this global definition of crossed modules yields naturally the notion of higher crossed modules. Namely, an $n\mhyphen$crossed$\mhyphen\p\mhyphen$algebra is a $\p\mhyphen$algebra structure on an $n\mhyphen$fold chain complex $C_{\bullet,...,\bullet}$ concentrated in degrees $\epsilon_1,...,\epsilon_n$ where $\epsilon_i \in \{0,1\}$. We prove that this notion of higher crossed modules is equivalent to the one of Ellis \cite{ELLIS1988277} for binary operads. Morevover, thanks to the monoidal totalization functor $Tot^\oplus: Ch(Ch(\A))\rightarrow Ch(\A)$, and multiple iterations of it, with $\A$ a bicomplete closed symmetric abelian category $\A$, we are able to map any $n\mhyphen$crossed$\mhyphen\p\mhyphen$algebra to a differential graded $\p\mhyphen$algebra concentrated in degrees $0$ to $n$.

\section*{Notations}In this paper, we fix a bicomplete abelian category $(\mathcal{A},\otimes,k)$ endowed with a closed symmetric monoidal structure, with unit of the tensor product denoted $k$. The biproduct is denoted $\oplus$, and it commutes with the tensor product on both variables as $\A$ is closed. The internal object of morphisms from $X$ to $Y$ is denoted $\underline{\mathcal{H}om}_\A(X,Y)$. If nothing is specified, we denote in the same way the enriched over itself category $\A$ or the underlying (set enriched) category $\A$ obtained by $Hom_\A(X,Y)\simeq Hom_\A(k,\Hom(X,Y))$. Chain complexes are considered over $\A$ and homologically $\mathbb{Z}\mhyphen$graded if nothing else is specified. This category of chain complexes is denoted $Ch(\A)$ and the full monoidal subcategory of non negatively graded chain complexes is denoted $Ch_+(\A)$. Given a symmetric operad $\mathcal{P}$ of $\mathcal{A}$, the object of $\mathcal{A}$ of $n\mhyphen$ary operations of $\mathcal{P}$ is denoted $\mathcal{P}_n$ and is endowed with a right action of the symmetric group $\Sigma_n$. In this paper we fix $\p$, a reduced ($\mathcal{P}_0 = 0 \text{ and } \mathcal{P}_1=k$ the unit of the tensor product) symmetric operad of $\A$. The operad $\p$ is said to be an $\A\mhyphen$operad. In particular the category $\mathcal{P}\mhyphen alg$ of algebras over this operad is bicomplete, has a zero object hence kernels. The operad $\p$ also defines a Schur endofunctor of $\A$ \cite{fresse2007modules, loday2012algebraic} denoted $\p()$. The operadic structure map of $\p$ induces a monad structure on $\p()$ and the multiplication of this monad $\p^2()\rightarrow\p()$ is denoted $\mu^{\p}$. In this paper we will make no difference between $\p$ and the Schur functor associated to it. For example, any object $x$ of $\A$ gives rise to the free $\p\mhyphen$algebra $\p(x)$ with structure map $\p^2(x)\rightarrow\p(x)$ induced by the structure map of $\p$. This free construction is left adjoint to the forgetful functor $U:\p\mhyphen\text{alg}\rightarrow \A$ from the category of $\p\mhyphen$algebras to the category $\A$. The theory of algebras over a symmetric operad $\mathcal{P}$ is particularly well behaved in the context of a symmetric \textit{closed} monoidal category because of the existence of the endomorphism operad $End(x)$ for every object $x \in \A$, which is given by $End(x)_n:=\Hom(x^{\otimes n},x)$. This operad allows to represent the contravariant functor that sends an operad $\p$ to the set of $\p\mhyphen$algebra structures on $x$ : \[\{\p\mhyphen\text{algebra structures on }x\}\simeq Hom_{Op}(\mathcal{P},End(x))\]
\section{Crossed Modules revisited}
This section is devoted to adapt the work of Janelidze \cite{janelidze2003internal} in the context where the category $\mathcal{C}$ of interest is the category of algebras over a reduced operad $\mathcal{P}$. His work takes place in the more general context of semi-abelian categories, and it is folklore that the category of algebras over a reduced operad, in modules over a commutative ring $R$, is so. Note that this context makes use of coproducts which are not explicit for a given operad. We then construct an equivalence of categories between crossed modules of $\p\mhyphen$algebras and internal category to $\p\mhyphen$algebras, adapting the method of Janelidze in this context.
%is not so appropriate when working with algebras over an operad, especially because it involves coproducts, which are not easy to deal with for algebras. The originality of this section is very limited and main ideas come from the work of Janelidze. Especially, we construct the equivalence of categories between crossed modules of algebras over an operad and internal categories to the category of such algebras in the same spirit as in \cite{janelidze2003internal}.

\subsection{Internal action}\label{section1.1}

\begin{definition}\label{actinterne}
Given two $\mathcal{P}\mhyphen$algebras $B$ and $X$, an \textit{internal action of $B$ on $X$} is the data of a morphism $\rho :\mathcal{P}\circ (B|X) \rightarrow X$ of the category $\A$, where :\[\mathcal{P}\circ (B|X) := \underset{n \in \mathbb{N}}\bigoplus \p_n\otimes_{\Sigma_n}((B \oplus X)^{\otimes n}/B^{\otimes n})\] This morphism must satisfy two properties : \begin{enumerate}
    \item(compatibility with the structure map of $X$) the triangle is commutative : \[\begin{tikzcd}
\mathcal{P}(X) \arrow[r, hook] \arrow[rd, "\mu^X"'] & \mathcal{P} \circ (B|X) \arrow[d, "\rho"] \\
                                                    & X                                        
\end{tikzcd}\]  Here the horizontal arrow is induced by $X \hookrightarrow X \oplus B$ and $\mu^X$ is the structure map of $X$. 

\item (associativity of the action) The following square is commutative : \[\begin{tikzcd}
\mathcal{P}\circ(\mathcal{P}(B)|\mathcal{P}\circ(B|X))  \arrow[d, "\mu^\mathcal{P}"'] \arrow[r] & \mathcal{P}\circ(B |X) \arrow[d, "\rho"] \\
\mathcal{P}\circ(B|X) \arrow[r, "\rho"']                                                                                                                                                 & X                                       
\end{tikzcd}\]Here the upper horizontal arrow is given by the restricion of $\mathcal{P}(\begin{pmatrix}\rho & 0 \\ 0 & \mu^B\end{pmatrix})$ to the summand $\mathcal{P}\circ(\mathcal{P}(B)|\mathcal{P}\circ(B|X))$ and $\mu^\mathcal{P}$ is the structure map of the operad $\mathcal{P}$.

\end{enumerate}Such an action is denoted $ B\overset{\rho}{\curvearrowright}X$.
\end{definition}
In the sequel of this paper, when dealing with an action without further details, this must be understood as an \textit{internal} action.
\begin{remarque}\label{actnat}
\begin{enumerate}
     \item $B\mhyphen \mathcal{P}\mhyphen $modules in the sense of \cite{fresse2007modules} are a particular case of such an action where the $\mathcal{P}\mhyphen $algebra $B$ acts on the trivial (any non identity operation vanishes) $\mathcal{P}\mhyphen $algebra $X$. In particular, an internal action of $B$ on $X$ gives rise to a $B\mhyphen\p\mhyphen$module structure on $X$ by forgetting the $\p\mhyphen$algebra structure on $X$.
    \item A $\mathcal{P}\mhyphen $algebra $B$ always acts on itself via the composition : \[\mathcal{P}\circ (B |B) \hookrightarrow \mathcal{P}(B \oplus B) \overset{\mathcal{P}(\mathrm{Id}\oplus \mathrm{Id})}{\longrightarrow} \mathcal{P}(B) \rightarrow B\]This action is denoted $\rho^B$.
    \item If $I \subset B$ is a $\p\mhyphen$ideal of $B$, then this natural action of $B$ on itself restricts to an action of $B$ on $I$. 
    \item Given an action $\mathcal{P}\circ(E|X) \rightarrow X$ of $E$ on $X$, and any morphism of $\p\mhyphen$algebras $f:B \rightarrow E$, the following composition $\mathcal{P}\circ(B|X)\overset{\mathcal{P}\circ(f|\mathrm{Id_X})}\longrightarrow \mathcal{P}\circ(E|X)\rightarrow X$ is a well defined action of $B$ on $X$.
\end{enumerate}
\end{remarque}
\begin{definition}
    A morphism between two actions $B \overset{\rho}{\curvearrowright}X$ and $B' \overset{\rho'}{\curvearrowright}X'$ is the data of a pair of $\mathcal{P}\mhyphen $algebra morphisms $(g: B \rightarrow B', f: X \rightarrow X')$ such that the following square is commutative : \[\begin{tikzcd}
\mathcal{P}\circ (B|X) \arrow[r, "\mathcal{P}\circ(g|f)"] \arrow[d, "\rho"'] & \mathcal{P}\circ (B'|X') \arrow[d, "\rho'"] \\
X \arrow[r, "f"']                                                            & X'                                         
\end{tikzcd}\]\end{definition}
\begin{definition}

The category $\mathrm{Act}(\p\mhyphen alg)$ is the category with objects the triples $(B,X,B \overset{\rho}{\curvearrowright}X)$. The morphisms of this category between two such triples $(B,X,B \overset{\rho}{\curvearrowright}X)$ and $(B',X',B \overset{\rho'}{\curvearrowright}X)$ are the morphisms of the previous definition.
    
\end{definition}
\begin{remarque}
    It is not difficult to derive from Remark \ref{actnat} that the assignement $B \mapsto \mathrm{Act}(B,X)$, which maps a $\p\mhyphen$algebra $B$ to the set of actions of $B$ on $X$, is a contravariant functor from the category of $\p\mhyphen$algebras to the category of sets. A result of \cite{GarciaVienne} states that, as soon as the category $\A$ is the category of $R\mhyphen$modules over an infinte field $R$ and $\p$ is generated by arity $2$ operations, this funtor is representable if and only if $\p\simeq\mathfrak{L}ie$. 
\end{remarque}
\begin{lemma}\label{act}
A split short exact sequence of $\mathcal{P}\mhyphen $algebras : \[\begin{tikzcd}
0 \arrow[r] & X \arrow[r, "i", hook] & E \arrow[r, "p", two heads] & B \arrow[r] \arrow[l, "s", bend left] & 0
\end{tikzcd}\]induces a natural internal action of $B$ on $X$.
    
\end{lemma}

\begin{demo}
The $\p\mhyphen$algebra morphism $i:X \rightarrow E$ is the kernel of $p$, so $X$ is a $\p\mhyphen$ideal of $E$, hence by $(3)$ of Remark \ref{actnat}, the natural action of $E$ on itself restricts to an action of $E$ on $X$. Moreover, as $s:B\rightarrow E$ is a $\p\mhyphen$algebra morphism, by $(4)$ of Remark \ref{actnat}, we get a well defined action of $B$ on $X$.

\end{demo}
\begin{lemma}\label{structurePalg}

An action of $B$ on $X$, $\rho : \mathcal{P} \circ (B | X) \rightarrow X$, induces a $\mathcal{P}\mhyphen $algebra structure on $X \oplus B$ by the formula :  \[\mathcal{P}(X \oplus B)=\mathcal{P}\circ(B|X) \oplus \mathcal{P}(B) \overset{\begin{pmatrix} \rho & 0 \\0 & \mu^B\end{pmatrix}}{\longrightarrow} X \oplus B\]

\end{lemma}

\begin{demo}
Let us denote this map by $\mu^{X \oplus B}$. We have to prove that $\mu^{X \oplus B}$ is associative, that is the following square commutes : \[\begin{tikzcd}
\mathcal{P}^2(X \oplus B) \arrow[r, "\mathcal{P}(\mu^{X \oplus B})"] \arrow[d, "\mu^\mathcal{P}"'] & \mathcal{P}(X \oplus B) \arrow[d, "\mu^{X \oplus B}"] \\
\mathcal{P}(X \oplus B) \arrow[r, "\mu^{X \oplus B}"']                                             & X \oplus B                                           
\end{tikzcd}\]
We know that $\mathcal{P}^2(X \oplus B)=\mathcal{P}\circ(\mathcal{P}(B)|\mathcal{P}\circ(B|X))\oplus \mathcal{P}^2(B)$, so we only have to check the commutativity of the previous square on each summand of $\mathcal{P}^2(X \oplus B)$. On the summand $\mathcal{P}^2(B)$, it is commutative because $B$ is a $\mathcal{P}\mhyphen$algebra, and on the other summand it is a direct consequence of the associativity of the action. 
\end{demo}
\begin{lemma}\label{secs}
    An action of $B$ on $X$ induces a split short exact sequence of $\mathcal{P}\mhyphen $algebras : \[\begin{tikzcd}
0 \arrow[r] & X \arrow[r, "i_1", hook] & X \oplus B \arrow[r, "p", two heads] & B \arrow[r] \arrow[l, "i_2", bend left] & 0
\end{tikzcd}\]Here the $\p\mhyphen$algebra structure on $X\oplus B$ is the one of Lemma \ref{structurePalg}.
\end{lemma}
\begin{demo}
    The compatibility of the action with the $\mathcal{P}\mhyphen $algebra structure on $X$ implies that $i_1$ is a $\mathcal{P}\mhyphen $algebra morphism. The definition of the $\mathcal{P}\mhyphen $algebra structure on $X \oplus B$ implies that the projection $p:X \oplus B \rightarrow B$ and $i_2$ are $\mathcal{P}\mhyphen $algebra morphisms. Exactness is obvious as the forgetful functor $U:\p\mhyphen alg\rightarrow A$ reflects exactness.
\end{demo}
\begin{lemma}
    A morphism of split short exact sequence of $\mathcal{P}\mhyphen$algebras induces a morphism between the corresponding actions given by Lemma \ref{act}.
\end{lemma}
\begin{demo}
    This is basic computation and it is left to the reader.
\end{demo}\newline

We denote by $\mathrm{SplitEpi}(\mathcal{P}\mhyphen alg)$ the category of split short exact sequences of $\p\mhyphen$algebras and its morphisms.
\begin{theorem}\label{actsplit}
    The constructions of Lemmas \ref{act} and \ref{secs} induce an equivalence of categories: \[\mathrm{Act}(\mathcal{P}\mhyphen alg) \simeq \mathrm{SplitEpi}(\mathcal{P}\mhyphen alg)\]
\end{theorem}

\begin{demo}
Previous constructions give both functors $\mathrm{Act}(\mathcal{P}\mhyphen alg) \rightarrow \mathrm{SplitEpi}(\mathcal{P}\mhyphen alg)$ and $\mathrm{SplitEpi}(\mathcal{P}\mhyphen alg) \rightarrow \mathrm{Act}(\mathcal{P}\mhyphen alg)$, and only the equivalence part of the theorem remains to be proved. On one side, given a split short exact sequence \[\begin{tikzcd}
0 \arrow[r] & X \arrow[r, hook] & E \arrow[r, "s", two heads] & B \arrow[r] \arrow[l, "i", bend left] & 0
\end{tikzcd},\] as $\A$ is an abelian category and $s$ is a split epimorphism, we have a natural isomorphism of split short exact sequences in $\A$ : \[\begin{tikzcd}
    0 \arrow[r] \arrow[d] & X \arrow[r, hook] \arrow[d, "\mathrm{Id}"'] & E \arrow[r, two heads, "s"] \arrow[d, "\simeq"]                    & B \arrow[r] \arrow[d, "\mathrm{Id}"] \arrow[l, "i"', bend left] & 0 \arrow[d] \\
    0 \arrow[r]           & X \arrow[r, "i_1", hook]               & X \oplus B \arrow[r, "p_2"] & B \arrow[r] \arrow[l, bend left, "i_2"]                     & 0          
    \end{tikzcd}\]Here $p_2$ and $i_2$ denote the canonical projection and inclusion respectively. Thus there is a natural $\p\mhyphen$algebra struture on $X \oplus B$, transfered from the one on $E$, such that this isomorphism of split short exact sequences is promoted to an isomorphism of split short exact sequences of $\p\mhyphen$algebras.

The action of $B$ on $X$ is given by the universal property of the $\mathcal{P}\mhyphen$ideal $X$: \[\begin{tikzcd}
\mathcal{P}\circ(B|X) \arrow[r, hook] \arrow[d, "\exists! \rho"', dotted] & \mathcal{P}(X\oplus B) \arrow[d, two heads, "\mu^{X \oplus B}"'] \arrow[r, "\mathcal{P}(p_2)"] & \mathcal{P}(B) \arrow[d, "\mu^B"] \arrow[l, "\mathcal{P}(i_2)", bend left] \\
X \arrow[r, hook]                                                         & X \oplus B \arrow[r, "p", two heads]                                   & B \arrow[l, "i_2", bend left]                                             
\end{tikzcd}\]We deduce that $\mu^{X \oplus B}= \begin{pmatrix}
    \rho & 0 \\ 0 & \mu^B
\end{pmatrix}$, in particular the composition : \[\mathrm{SplitEpi}(\mathcal{P}\mhyphen alg) \rightarrow \mathrm{Act}(\mathcal{P}\mhyphen alg) \rightarrow \mathrm{SplitEpi}(\mathcal{P}\mhyphen alg)\]is naturally isomorphic to the identity. On the other side, it is clear that the composition $\mathrm{Act}(\mathcal{P}\mhyphen alg) \rightarrow \mathrm{SplitEpi}(\mathcal{P}\mhyphen alg) \rightarrow \mathrm{Act}(\mathcal{P}\mhyphen alg)$ is naturally isomorphic to the identity. 
    
\end{demo}

\subsection{Precrossed modules}

\begin{definition}
\begin{itemize}

   \item A \textit{precrossed module} is a pair $(B \overset{\rho}{\curvearrowright}X, d: X \rightarrow B)$ where $B$ acts on $X$ by $\rho$ and $d$ is a $B\mhyphen$equivariant morphism of $\p\mhyphen$algebras, for the natural action $\rho^B$ of Remark \ref{actnat} of $B$ on itself.

   \item A \textit{morphism of precrossed modules} is a morphism of actions commuting with the equivariant morphisms : \[\begin{tikzcd}
X \arrow[r, "d"] \arrow[d, "f"'] & B \arrow[d, "g"] \\
X' \arrow[r, "d'"']              & B'              
\end{tikzcd}\]
   \end{itemize}We denote by $\mathrm{pXMod}(\mathcal{P}\mhyphen alg)$ the category of precrossed module of $\mathcal{P}\mhyphen $algebras.
\end{definition}
\begin{remarque}
    Recall from Lemma \ref{structurePalg} that an action of $B$ on $X$ gives rise to a $\p\mhyphen$algebra structure on $X\oplus B$. We can check that in this context $d:X \rightarrow B$ is $B\mhyphen$equivariant if and only if $d \oplus \mathrm{Id} : X \oplus B \rightarrow B$ is a morphism of $\mathcal{P}\mhyphen $algebras.
\end{remarque}
In section \ref{section1.1}, we established an equivalence of categories between the category of internal actions and the category of split short exact sequences, which restricts to an equivalence between actions $B \overset{\rho}{\curvearrowright}X$ and the category of split short exact sequences with kernel $X$ and cokernel $B$. When the action is endowed with an equivariant morphism $d:X \rightarrow B$, we will see that it naturally leads to consider another morphism of $\mathcal{P}\mhyphen $algebras $d \oplus \mathrm{Id}:X \oplus B \rightarrow B$.

\begin{definition}\label{defreflgraph}
    A \textit{reflexive graph} of $\mathcal{P}\mhyphen $algebras is the data of a diagram of the following form in the category of $\mathcal{P}\mhyphen $algebras : \[\begin{tikzcd}
E \arrow[r, "s", shift left=2] \arrow[r, "t"', shift right=2] & B \arrow[l, "i" description]
\end{tikzcd}\]such that $ti=si=\mathrm{Id}_B$. A \textit{morphism of reflexive graphs} is a morphism of diagrams. We denote $\mathrm{ReflGraph}(\mathcal{P}\mhyphen alg)$ the category of reflexive graphs of $\mathcal{P}\mhyphen $algebras.
\end{definition}

Recall from the proof of Theorem \ref{actsplit} that there is an isomorphism of split short sequences of $\p\mhyphen$algebras: \[\begin{tikzcd}
0 \arrow[r] \arrow[d] & \mathrm{ker}(s) \arrow[r, hook] \arrow[d, "="'] & E \arrow[r, "s"] \arrow[d, "\simeq"]                    & B \arrow[r] \arrow[d, "\mathrm{Id}"] \arrow[l, "i"', bend left] & 0 \arrow[d] \\
0 \arrow[r]           & X \arrow[r, "i_1", hook]               & X \oplus B \arrow[r, "p_2"] & B \arrow[r] \arrow[l, bend left, "i_2"]                     & 0          
\end{tikzcd}\]Here $X$ denotes the kernel of $s$. We denote the composition $X \oplus B \overset{\simeq}\rightarrow E \overset{t}{\rightarrow} B $ by $d\oplus \mathrm{Id}$. This is a $\p\mhyphen$algebra morphism, hence $d:=t_{|\mathrm{ker}(s)}: X \rightarrow B$ is also a $\p\mhyphen$algebra morphism. \\
\newline Forgetting the morphism $t$ of a reflexive graph as in Definition\ref{defreflgraph}, gives rise to a split short exact sequence of $\p\mhyphen$algebras. Thus by Lemma \ref{act}, there is a natural internal action $B \overset{\rho}{\curvearrowright}X$.

\begin{lemma} Given a reflexive graph \begin{tikzcd}
X \oplus B \arrow[r, "0\oplus \mathrm{Id}", shift left=2] \arrow[r, "d\oplus \mathrm{Id}"', shift right=2] & B \arrow[l, "i_2" description]
\end{tikzcd} of $\mathcal{P}\mhyphen $algebras, the morphism $d: X \rightarrow B$ is equivariant with respect to the action of $B$. This provides a functor from the category of reflexive graphs to the category of precrossed modules.
\end{lemma}
\begin{demo}
    This is a direct consequence of the fact that $d$ and $d \oplus \mathrm{Id}$ are morphisms of $\mathcal{P}\mhyphen $algebras.
\end{demo}
\begin{lemma}
    Given a precrossed module $(B \overset{\rho}{\curvearrowright}X, d : X \rightarrow B)$, the morphism $d \oplus \mathrm{Id}: X \oplus B \rightarrow B$ is a morphism of $\mathcal{P}\mhyphen $algebras for the natural structure of $\p\mhyphen$algebra on $X \oplus B$ induced by Lemma \ref{structurePalg}.
\end{lemma}

\begin{demo} Recall from Lemma \ref{structurePalg} that the action of $B$ on $X$ induces a $\p\mhyphen$algebra structure on $X\oplus B$, that is a map $\mu^{X\oplus B}:\p(X\oplus B) \rightarrow X \oplus B$.
    To prove that $d\oplus \mathrm{Id}$ is a morphism of $\p\mhyphen$algebras, we have to prove that the following square is commutative :\[\begin{tikzcd}
\mathcal{P}(X \oplus B) \arrow[d,"\mu^{X\oplus B}"] \arrow[r, "\mathcal{P}(d \oplus \mathrm{Id})"] & \mathcal{P}(B) \arrow[d, "\mu^B"] \\
X \oplus B \arrow[r, "d \oplus \mathrm{Id}"']                                                                                            & B                                
\end{tikzcd}\] The decomposition $\mathcal{P}(X \oplus B)= \mathcal{P}\circ (B|X) \oplus \mathcal{P}(B)$ implies that it is enough to prove that the outer square of the following diagram commutes : \[\begin{tikzcd}
\mathcal{P}\circ (B |X) \arrow[dd, "\rho"'] \arrow[rd, "\mathcal{P}(B|d)" description] \arrow[rr, "\mathcal{P}(d \oplus \mathrm{Id})"] &                                                                                                & \mathcal{P}(B) \arrow[dd, "\mu^B"] \\
                                                                                                                              & \mathcal{P}\circ(B|B) \arrow[ru, "\mathcal{P}(\mathrm{Id} \oplus \mathrm{Id})" description] \arrow[rd, "\rho^B"] &                                    \\
X \arrow[rr, "d"']                                                                                                            &                                                                                                & B                                 
\end{tikzcd}\]

As $d$ is equivariant the bottom left part commutes, the right triangle also commutes by definition of the natural action of $B$ on itself. The upper triangle commutes thanks to the definition of the restriction.  
\end{demo}

Thus we deduce the following theorem.
\begin{theorem}\label{precross}
    The equivalence of categories of Theorem \ref{actsplit} lifts to an equivalence of categories $\mathrm{pXMod}(\mathcal{P}\mhyphen alg) \simeq \mathrm{ReflGraph}(\mathcal{P}\mhyphen alg)$ making the diagram commute: \[\begin{tikzcd}
\mathrm{Act}(\mathcal{P}\mhyphen alg) \arrow[r, "\simeq"]                        & \mathrm{SplitEpi}(\mathcal{P}\mhyphen alg)                  \\
\mathrm{pXMod}(\mathcal{P}\mhyphen alg) \arrow[u, "U"] \arrow[r, "\simeq"] & \mathrm{ReflGraph}(\mathcal{P}\mhyphen alg) \arrow[u, "U"']
\end{tikzcd}\]Here the left vertical arrow is given by forgetting the equivariant map and the right vertical one is given by forgetting the target arrow $t$.
\end{theorem}

\subsection{Crossed modules}\label{modulecroisé}
The main point of this section is to prove that on a reflexive graph \begin{tikzcd}
E \arrow[r, "s", shift left=2] \arrow[r, "t"', shift right=2] & B \arrow[l, "i" description]
\end{tikzcd}, there exists a composition $\circ : E \times_{s,t} E \rightarrow E$, turning it into an internal category to the category of $\mathcal{P}\mhyphen $algebras, if and only if $\forall n \geq 2$ the following diagram is commutative : \[\begin{tikzcd}
\mathcal{P}_n\otimes(\mathrm{ker}(s)\otimes \mathrm{ker}(t)\otimes E^{\otimes n-2}) \arrow[r] \arrow[rd, "0"'] & \mathcal{P}_n\underset{\Sigma_n}{\otimes}E^{\otimes n} \arrow[d, "\mu^E"] \\
                                                                                                                        & E                                                                        
\end{tikzcd}\]

\begin{definition}\label{Peiffer}
    A precrossed module $(B \overset{\rho}{\curvearrowright}X, d : X \rightarrow B)$ of $\mathcal{P}\mhyphen $algebras is a \textit{crossed module} if the $\mathcal{P}\mhyphen $algebra structure on $X$ is already determined by both the $B\mhyphen \mathcal{P}\mhyphen $module structure on $X$ and $d$. That is $\forall n,i \text{ such that } n \geq i \geq 2$, we have the following commutative diagram, called the \textit{Peiffer identity} : \begin{equation}\begin{tikzcd}
\mathcal{P}_n\otimes(X^{\otimes i}\otimes B^{\otimes n-i}) \arrow[r] \arrow[d, "\mathrm{Id}\otimes \mathrm{Id}\otimes d^{\otimes i-1}\otimes \mathrm{Id}"'] & \mathcal{P}\circ (B|X) \arrow[rd, "\rho"]  &   \\
\mathcal{P}_n\otimes(X\otimes B^{\otimes n-1} ) \arrow[rd]                                                                                           &                                            & X \\
                                                                                                                                                                                & \mathcal{P}\circ (B|X) \arrow[ru, "\rho"'] &  
\end{tikzcd}\label{eq:peiffer} \tag{2}\end{equation} We denote by $\mathrm{XMod}(\p\mhyphen alg)$, the full subcategory of precrossed modules, with crossed modules as objects.
\end{definition}

\begin{remarque}
    If $\mathcal{A}$ is a concrete category, then the commutativity of diagram (\ref{eq:peiffer}) is equivalent to the following condition : $\forall \gamma \in \mathcal{P}_n, n\geq 2$, we have the equality : \[\rho(\gamma;x_1,...,x_i,b_{i+1},...,b_n)=\rho(\gamma;x_1,dx_2,...,dx_i,b_{i+1},...,b_n)\]Here $n \geq i \geq 2$, $x_k$'s are in $X$ and $b_k$'s are in $B$.
\end{remarque}
\begin{definition}\label{cat^1}
    A reflexive graph of $\mathcal{P}\mhyphen $algebras \begin{tikzcd}
E \arrow[r, "s", shift left=2] \arrow[r, "t"', shift right=2] & B \arrow[l, "i" description]
\end{tikzcd} is said to be a \textit{ $\mathrm{Cat}^1\mhyphen \mathcal{P}\mhyphen $algebra} if $\forall n \geq 2$ the following diagramm commutes : \[\begin{tikzcd}
\mathcal{P}_n\otimes(\mathrm{ker}(s)\otimes \mathrm{ker}(t)\otimes E^{\otimes n-2}) \arrow[r] \arrow[rd, "0"'] & \mathcal{P}(E) \arrow[d, "\mu^E"] \\
                                                                                                                        & E                                
\end{tikzcd}\] This condition is called the $\mathrm{Cat}^1\mhyphen$condition. The full subcategory of reflexive graphs of $\mathcal{P}\mhyphen $algebras spanned by $\mathrm{Cat}^1\mhyphen \mathcal{P}\mhyphen $algebras is denoted $\mathrm{Cat}^1(\mathcal{P}\mhyphen alg)$.
\end{definition}
\begin{remarque}
\begin{itemize}
    \item Note that as $\Sigma_n$ acts on $E^{\otimes n}$, one can have equally asked the same equality with $\mathcal{P}_n\otimes(\mathrm{ker}(t)\otimes \mathrm{ker}(s)\otimes E^{\otimes n-2})$ or any other permutation of the factors of $\mathrm{ker}(s)\otimes \mathrm{ker}(t)\otimes E^{\otimes n-2}$, so there is no lack of symmetry.
    \item As $\mathrm{ker}(s)$ is a $\mathcal{P}\mhyphen $ideal of $E$, the commutativity of the previous diagram is equivalent to the commutativity of the following diagram : \[\begin{tikzcd}
\mathcal{P}_n\otimes(\mathrm{ker}(s)\otimes \mathrm{ker}(t)\otimes E^{\otimes n-2}) \arrow[rd, "0"'] \arrow[r] & \mathcal{P}(E) \arrow[d, "\rho"] \\
                                                                                                                      & \mathrm{ker}(s) \cap \mathrm{ker}(t)\subset E                                       
\end{tikzcd}\]
    \item If $s=t$, then this $\mathrm{Cat}^1\mhyphen$condition implies that the kernel of $s$ is a trivial $\p\mhyphen$algebra. In particular $E$ is a split abelian extension of $B$ by $X:=\mathrm{ker}(s)$.
    \end{itemize}
\end{remarque}
\begin{proposition}\label{peiffer=Cat^1}
The equivalence of categories between precrossed modules and reflexive graphs restricts to an equivalence of categories : \[\mathrm{XMod}(\p\mhyphen alg)\simeq \mathrm{Cat}^1(\mathcal{P}\mhyphen alg)\]
\end{proposition}
\begin{demo}

Given a crossed module $d:X \rightarrow B$, we want to prove that the reflexive graph we get : \begin{tikzcd}
X \oplus B \arrow[r, "s:=0 \oplus \mathrm{Id}", shift left=2] \arrow[r, "t:=d \oplus \mathrm{Id}"', shift right=2] & B \arrow[l, "i_2" description]
\end{tikzcd} is a $\mathrm{Cat}^1\mhyphen \mathcal{P}\mhyphen $algebra, that is, the restriction of the structural map $\mathcal{P}_n\otimes(X\otimes \mathrm{ker}(t)\otimes E^{\otimes n-2}) \rightarrow X$ is $0$. The map $\begin{tikzcd}
X \arrow[r, "\binom{\mathrm{Id}}{-d}"] & X \oplus B
\end{tikzcd}$ is a kernel of $t:=d\oplus \mathrm{Id} : X \oplus B \rightarrow B$, thus we get a well defined isomorphism : \[
\mathcal{P}_n\otimes(X\otimes X \otimes E^{\otimes n-2}) \overset{{\mathrm{Id}\otimes(\mathrm{Id}\otimes \binom{\mathrm{Id}}{-d}\otimes \mathrm{Id})}}{\longrightarrow}\mathcal{P}_n\otimes(X\otimes \mathrm{ker}(t) \otimes E^{\otimes n-2})\]So we are led to prove that the composition $\rho \circ (\mathrm{Id}\otimes(\mathrm{Id}\otimes \binom{\mathrm{Id}}{-d}\otimes \mathrm{Id})) :  \mathcal{P}_n\otimes(X\otimes X \otimes E^{\otimes n-2}) \rightarrow X$ is the $0$ map. But as the tensor product commutes with biproducts (as a left adjoint on both sides) and the composition is bilinear, the composition being $0$ is equivalent to $\rho=\rho\circ({\mathrm{Id}\otimes(\mathrm{Id}\otimes d\otimes \mathrm{Id})})$, or expressed in terms of diagrams, it is equivalent to the commutativity of the following pentagon : \[\begin{tikzcd}
\mathcal{P}_n\otimes(X^{\otimes 2}\otimes E^{\otimes n-2}) \arrow[r] \arrow[d, "1\otimes (1\otimes d \otimes 1)"'] & \mathcal{P}\circ (B|X) \arrow[rd, "\rho"]  &   \\
\mathcal{P}_n\otimes(X\otimes E^{\otimes n-1} ) \arrow[rd]                                                                                           &                                            & X \\
                                                                                                                                                                                & \mathcal{P}\circ (B|X) \arrow[ru, "\rho"'] &  
\end{tikzcd}\]
The commutativity follows from the Peiffer identity fulfilled by the crossed module\newline $X\overset{d}{\rightarrow} B$.

On the other hand, given a $\mathrm{Cat}^1\mhyphen \mathcal{P}\mhyphen $algebra \begin{tikzcd}
E \arrow[r, "s", shift left=2] \arrow[r, "t"', shift right=2] & B \arrow[l, "i" description]
\end{tikzcd}, to prove that the precrossed module $d:=t_{|\mathrm{ker}(s)}:\mathrm{ker}(s)\rightarrow B$ satisfies the Peiffer identity (\ref{eq:peiffer}) is very similar to the previous part of the proof, using induction on $i$, and it is left to the reader.
\end{demo}
\begin{theorem}\label{cat^1=cat}
    The identity on objects and morphisms is a functor which is an isomorphism of categories : \[\mathrm{Cat}^1(\mathcal{P}\mhyphen alg)\simeq \mathrm{Cat}(\mathcal{P}\mhyphen alg)\]
\end{theorem}
\begin{demo}
    For the sake of simplicity, let us assume that $\mathcal{A}$ is the abelian category of $k\mhyphen $modules, or chain complexes of $k\mhyphen $modules, endowed with its closed monoidal structure provided by the usual tensor product of chain complexes over the commutative ring $k$. This viewpoint allows us to pick elements instead of writing down big diagrams which are not so enlightening. We will use nothing about $k\mhyphen $modules except that the tensor product commutes with biproducts, which is assumed in our abelian category $\mathcal{A}$, and that in abelian categories, an internal reflexive graph can be uniquely promoted to an internal category \cite{bourn1990another}. In particular, the free-forgetful adjunction \begin{tikzcd}
            \mathcal{P}\mhyphen alg\arrow[r, shift left=1ex, "U"{name=G}] & \mathcal{A}\arrow[l, shift left=.5ex, "\mathcal{P}(\mhyphen )"{name=F}]
            \arrow[phantom, from=F, to=G, , "\scriptscriptstyle\boldsymbol{\top}"]
        \end{tikzcd}
implies that on a reflexive graph internal to $\mathcal{P}\mhyphen $algebras, there exists at most one composition which turns it into an internal category to the category of $\mathcal{P}\mhyphen $algebras.
First, let us prove that internal categories give rise to $\mathrm{Cat}^1\mhyphen \mathcal{P}\mhyphen $algebras. Let \begin{tikzcd}
E \arrow[r, "s", shift left=2] \arrow[r, "t"', shift right=2] & B \arrow[l, "i" description]
\end{tikzcd} be an internal category to $\mathcal{P}\mhyphen $algebras and consider $\gamma \in \mathcal{P}_n, f_2 \in \mathrm{ker}(t),g_1 \in \mathrm{ker}(s)$, $f_1:=\mathrm{Id}_0=0$ and $g_2:=\mathrm{Id}_0=0$. In particular $g_2 \circ f_2$ and $g_1 \circ f_1$ are well defined and we can compute, for $\bar{a}\in E^{\otimes n-2}$ : \begin{align*}
        \gamma(g_1,f_2,\bar{a}) & = \gamma(g_1 \circ f_1, g_2 \circ f_2, \bar{a}\circ \mathrm{Id}) \\ & = \gamma(g_1,g_2, \bar{a}) \circ \gamma(f_1,f_2, \mathrm{Id}) \\ & = \gamma(g_1,0,\bar{a}) \circ \gamma(0,f_2,\mathrm{Id}) \\ &= 0 \circ 0 \\ & = 0
\end{align*} The transition from the first line to the second one is allowed because the composition is a morphism of $\mathcal{P}\mhyphen $algebras. This proves that internal categories give rise to $\mathrm{Cat}^1\mhyphen\p\mhyphen$algebras. \newline On the other side, we have to prove that if \begin{tikzcd}
E \arrow[r, "s", shift left=2] \arrow[r, "t"', shift right=2] & B \arrow[l, "i" description]
\end{tikzcd} is a $\mathrm{Cat}^1\mhyphen \mathcal{P}\mhyphen $algebra, then the only possible composition on it, given by $g \circ f = f + g - isg$, is a morphism of $\mathcal{P}\mhyphen $algebras. Given $\gamma \in \mathcal{P}(n), n \geq 2$, and pairs of composable morphisms $(g_i,f_i)_{1 \leq i \leq n} \in E \times_{B} E$, we have to prove the following equality: \[\gamma(g_1 \circ f_1, ...,g_n \circ f_n)=\gamma(g_1,...,g_n) \circ \gamma(f_1,...,f_n)\]
First, we assume that each $f_k$ is in the kernel of $t$, so $g_k$'s are in the kernel of $s$. Thanks to the formula for the composition, we know that $g_k \circ f_k=g_k + f_k$, so it allows us to compute : \begin{align*}
    \gamma(g_1\circ f_1, ... , g_n \circ f_n) & = \gamma(g_1 + f_1,..., g_n + f_n) \\ & = \gamma(g_1,...,g_n) + \gamma(f_1,...,f_n) \\ & = \gamma(g_1,...,g_n) \circ \gamma(f_1,...,f_n)
\end{align*} Indeed, transition from the first to the second line is because $\gamma$ is $n\mhyphen $linear, thus by expanding the product, terms $\gamma(a_1,...,a_n)$ where at least one $f_k$ and one $g_j$ appear are $0$ because of the $\mathrm{Cat}^1\mhyphen \mathcal{P}\mhyphen $algebra condition : $\gamma(\mathrm{ker}(s),\mathrm{ker}(t),-)=0$. And transition from second to last line is due to the fact that $s$ and $t$ are morphisms of $\mathcal{P}\mhyphen $algebras, in particular $s(\gamma(g_1,...,g_n))=t(\gamma(f_1,...,f_n))=0$.\newline Now assume that some composable pairs $(g_i,f_i)$ are of the form $(i(x_k),i(x_k))$. Without loss of generality, up to an action of $\Sigma_n$, one can assume that the $k$ first composable pairs are of this form. If $k=n$, then :\begin{align*}
    \gamma(i(x_1)\circ i(x_1),...,i(x_n)\circ i(x_n)) & = \gamma(i(x_1),...,i(x_n)) \\ & = \gamma \circ i (x_1,...,x_n) \\ & = i(\gamma(x_1,...,x_n)) \\ & = i(\gamma(x_1,...,x_n)) \\ & = i(\gamma(x_1,...,x_n)) \circ i(\gamma(x_1,...,x_n)) \\ & = \gamma(i(x_1),...,i(x_n)) \circ \gamma(i(x_1),...,i(x_n))
\end{align*}Now let us assume that $k<n$ and $\forall j \geq k+1 : sg_j=tf_j=0$, that allows us to compute  \begingroup  \begin{align*}
    & \gamma(i(x_1)\circ i(x_1),...,i(x_k)\circ i(x_k),g_{k+1}\circ f_{k+1},...,g_n \circ f_n) \\ & = \gamma(i(x_1),...,i(x_k),g_{k+1}\circ f_{k+1},...,g_n \circ f_n) \\ &= \gamma(i(x_1),...,i(x_k),g_{k+1}+ f_{k+1},...,g_n + f_n)\\ & = \gamma(i(x_1),...,i(x_k),g_{k+1},...,g_n) + \gamma(i(x_1),...,i(x_k),f_{k+1},...,f_n) \\ & = \gamma(i(x_1),...,i(x_k),g_{k+1},...,g_n) \circ \gamma(i(x_1),...,i(x_k),f_{k+1},...,f_n)
\end{align*}\endgroup Equality between the third line and the fourth one comes from the multilinearity of $\gamma$ (a.k.a. tensor product commutes with coproducts) and the $\mathrm{Cat}^1\mhyphen \mathcal{P}\mhyphen $algebra condition. Equality between the fourth line and the fifth one comes from the fact that $k<n$, so the source of $\gamma(\mathrm{Id}_{x_1},...,\mathrm{Id}_{x_k},g_{k+1},...,g_n)$ and the target of $\gamma(\mathrm{Id}_{x_1},...,\mathrm{Id}_{x_k},f_{k+1},...,f_n)$ are $0$, thus the sum of these two terms is precisely their composition. Now let us prove by induction on $l$ the property $P_l$ : $\gamma(g_1\circ f_1,...,g_n \circ f_n)=\gamma(g_1,...,g_n) \circ \gamma(f_1,...,f_n)$, where $f_i,g_i$ are composable and :\[l:=\text{Card}(\llbracket 1,n\rrbracket\textbackslash \{ i\in \llbracket 1,n\rrbracket|(sg_i=tf_i=0) \vee ((sg_i=tf_i\neq 0 )\Rightarrow sg_i=tf_i=i(x_i))\})\] The number $l$ counts the number of composable pairs which are not of the form $(i(x_i), i(x_i))$ nor $(g_i,f_i)$ with $sg_i=tf_i=0$. The previous work is the initialisation step $P_0$. Assuming $P_l$ is true, we want to prove that $P_{l+1}$ is also true. Without loss of generality, we can assume that $(g_1,f_1)$ is neither a couple of identity nor morphisms with source and target $0$. One can write $ g_1=g_1-i(sg_1) + i(sg_1)$, with $g_1-i(sg_1)$ having $0$ as source. One can also write $f_1=f_1-i(tf_1) + i(tf_1)$ with $f_1-i(f_1)$ having $0$ as target. Let us write $\bar{g}:=(g_2,...,g_n), \bar{f}:=(f_1,...,f_n), \bar{g}\circ \bar{f}:=(g_2 \circ f_2,...,g_n \circ f_n)$ and $x_1:=sg_1=tf_1$, and the computation gives : \begin{align*}
    \gamma(g_1 \circ f_1,\bar{g}\circ \bar{f}) & = \gamma ((g_1-i(x_1)+i(x_1))\circ(f_1-i(x_1) + i(x_1)),\bar{g}\circ \bar{f}) \\ & = \gamma((g_1-i(x_1))\circ (f_1-i(x_1)) + (i(x_1)\circ i(x_1)), \bar{g}\circ \bar{f}) \\ & = \gamma((g_1-i(x_1)) \circ (f_1-i(x_1)), \bar{g}\circ \bar{f}) + \gamma(i(x_1)\circ i(x_1), \bar{g}\circ\bar{f}) \\ & = \gamma(g_1-i(x_1),\bar{g})\circ \gamma(f_1-i(x_1),\bar{f}) + \gamma(i(x_1),\bar{g}) \circ \gamma(i(x_1),\bar{f}) \\ & = \gamma(g_1,\bar{g})\circ \gamma(f_1,\bar{f})
\end{align*}Induction hypothesis is used on both factors on the transition from third to fourth line. Transition from fourth to fifth is by linearity of the composition as well as transition from first to second line. This concludes the proof.
\end{demo}
\subsection{Conclusion}
To conclude, we have proven that various categories, defined internally to the category of $\mathcal{P}\mhyphen $algebras, are isomorphic. It is compiled in the following commutative diagram, where vertical arrows are forgetful functors : \[\begin{tikzcd}
\mathrm{Act}(\mathcal{P}\mhyphen alg) \arrow[r, "\simeq"]                   & \mathrm{SplitEpi}(\mathcal{P}\mhyphen alg)            \\
\mathrm{pXMod}(\mathcal{P}\mhyphen alg) \arrow[u] \arrow[r, "\simeq"] & \mathrm{ReflGraph}(\mathcal{P}\mhyphen alg) \arrow[u] \\
\mathrm{XMod}(\p\mhyphen alg) \arrow[u] \arrow[r, "\simeq"]            & \mathrm{Cat}(\mathcal{P}\mhyphen alg) \arrow[u]    
\end{tikzcd}\]
Here, vertical arrows are forgetful functors and, moreover, the bottom ones are fully faithful embeddings. Note that the category $\mathrm{Cat}^1(\mathcal{P}\mhyphen alg)$ is isomorphic to the category $\mathrm{Cat}(\mathcal{P}\mhyphen alg)$ of internal categories to the category of $\mathcal{P}\mhyphen $algebras, so throughout the paper, we will use equally these two notions as internal categories have a more geometric interpretation and the $\mathrm{Cat}^1$ description is great to compute with.

\section{Crossed $\mathcal{P}\mhyphen$algebras}
Through the equivalence $\mathrm{XMod}(\p\mhyphen alg) \simeq \mathrm{Cat}(\p\mhyphen alg)$ and in the light of Loday's work \cite{loday1982spaces} on homotopy $n\mhyphen$types, there is an evident generalization for higher crossed modules of $\p\mhyphen$algebras, namely $n\mhyphen$fold categories internal to the category of $\p\mhyphen$algebras, which we denote $\mathrm{Cat}^n(\p\mhyphen alg)$. Then one can ask what can play the role of $n\mhyphen$crossed modules, denoted $\mathrm{X}^nMod(\p\mhyphen alg)$, in the aim of having an equivalence of categories $\mathrm{X}^nMod(\p\mhyphen alg) \simeq \mathrm{Cat}^n(\p\mhyphen alg)$. A first answer is to consider crossed modules internal to $(n-1)\mhyphen$crossed modules of $\p\mhyphen$algebras. For example, a $2\mhyphen$crossed module should be the data of two crossed modules of $\p\mhyphen$algebras, say $L\rightarrow M$ and $N\rightarrow P$, a morphism between them, that is, a commutative square of $\mathcal{P}\mhyphen$algebras :\[\begin{tikzcd}
L \arrow[d] \arrow[r] & N \arrow[d] \\
M \arrow[r]           & P          
\end{tikzcd}\] This should be endowed with an action of $N\rightarrow P$ on $L\rightarrow M$ such that together with the horizontal arrows, they form an equivariant map. This morphism of crossed modules must also satisfy a Peiffer condition. The work for some particular cases of algebras has been done by Ellis in \cite{ELLIS1988277}, for example in the cases of Lie/associative/commutative algebras. There are a few issues with his approach. First, equations for higher crossed modules are quite involved. Then it is not obvious how to develop a notion of homotopy theory for these objects, and if possible to relate it with the homotopy theory of $\p\mhyphen$algebras, namely differential graded $\p\mhyphen$algebras. One of the goals of this paper is to fill this gap. The starting point is the following. In a current work \cite{lrw}, Leray-Rivière-Wagemann are considering another notion of crossed modules of $\p\mhyphen$algebras : it is the data of a $\p\mhyphen$algebra structure on a chain complex $V:=...0\rightarrow V_1 \rightarrow V_0\rightarrow 0 ... $ concentrated in homological degrees $0$ and $1$. The advantage of this approach is two fold : it is clearly linked with the homotopy theory of $\p\mhyphen$algebras and it is also a very concise definition. But the difficulty has been moved : we have to show that we recover the classical definition of a crossed module. It is the goal of the first subsection by enlightening the "derived operation" procedure of Kosmann-Schwarzbach\cite{Kosmann_Schwarzbach_2004}, or codescent phenomenon, which allows us to induce a $\p\mhyphen$algebra structure on $V_1$ such that $V_1 \overset{d}{\rightarrow} V_0$ is a crossed module of $\mathcal{P}\mhyphen$algebras. This global point of view, as opposed to the local one from Ellis or Janelidze, allows a very concise formalism for higher crossed modules : they should be $\p\mhyphen$algebra structures on $n\mhyphen$fold chain complexes $C_{\bullet,...,\bullet}$ concentrated in degrees $(\epsilon_1,...,\epsilon_n)$, where $\epsilon_i \in \{0,1\}$. The category of these objects, called $n\mhyphen$crossed modules of $\p\mhyphen$algebras, with $\p\mhyphen$algebra morphisms as morphisms, is denoted $\mathrm{X}^n\mathrm{Mod}(\p\mhyphen alg)$. Using the codescent phenomenon, we recover $\p\mhyphen$algebra structures on each $C_{\epsilon_1,...,\epsilon_n}$. Iterations of the totalisation functor $\mathrm{Tot}: Ch_+(Ch_+(\A))\rightarrow Ch_+(\A)$, which is monoidal, allow us to associate to any $n\mhyphen$crossed module of $\p\mhyphen$algebras a differential graded $\p\mhyphen$algebra concentrated in degrees $0$ to $n$. Hence $\mathrm{Tot}$ provides a bridge between our notion of $n\mhyphen$crossed modules and the homotopy theory of $\p\mhyphen$algebras.

\subsection{Iterated chain complexes}
We refer to Kelly \cite{kelly1982basic} for a global account on enriched category theory. This subsection is devoted to recall some basic and folklore facts about chain complexes in a closed symmetric monoidal abelian category in order to treat $Ch(\A)$ and $n\mhyphen$fold chain complexes $Ch(...(Ch(\A)))$ (where "$Ch$" appears $n\mhyphen$times) with the same formalism. We denote the category of $n\mhyphen$fold chain complexes on $\A$ by $nCh(\mathcal{A})$.The key point is that if $\A$ is a closed symmetric monoidal category, then so are $Ch(\A)$ and $Ch_+(\A)$. Recall that in an abelian category, the set of morphisms between two objects is an abelian group. In the closed symmetric monoidal category $\A$, as every object is an internal abelian group object, so is the (internal) object $\Hom(x,y)$, that is $\Hom(x,y)\in Ab(\A)$. Indeed the abelian group structure is given by the following data: \begin{itemize}
     \item  The group law is $m:\Hom(x,y)\oplus \Hom(x,y) \overset{\mathrm{Id}\oplus \mathrm{Id}}{\rightarrow} \Hom(x,y)$
\item The zero element is given by the unique map $0\overset{\exists !}{\rightarrow} \Hom(x,y)$
\item The inversion $Inv:\Hom(x,y) \rightarrow \Hom(x,y)$ is given by $-\mathrm{Id}_{\Hom(x,y)}$
 \end{itemize} Moreover $(k,\oplus,\otimes,k\overset{\mathrm{Id}}{\rightarrow}k)$ is a unital ring object internal to $\A$ and every object of $\A$ carries a natural $k\mhyphen$module structure, thus $\A\simeq \mathrm{Mod}_k$.
\newline Here, a graded object of A is $\mathbb{Z}$-graded. It consists of a collection of objects $(V_n)_{n \in \mathbb{Z}}$ of $\A$ and we denote it $V_\bullet$ or simply $V$. Equivalently, it is a functor $V:\mathbb{Z}\rightarrow \A$ from the discrete category $\mathbb{Z}$ to the category $\A$. We denote by $\mathrm{gr}\A$ the category of graded objects and natural transformations. This category has automorphisms induced by shifting graded objects. We denote $V_{\bullet -1}$ the object $V$ shifted by $-1$, that is $(V_{\bullet -1})_n=V_{n-1}$. The closed symmetric monoidal structure of $\mathrm{gr}\A$ is constructed as follows. For $r \in \mathbb{Z}$, the object of \textit{morphisms of degree r} between two graded objects $A$ and $B$, denoted by $\underline{Hom}(A,B)_r$, is the product $\underset{n\in \mathbb{Z}}{\Pi}\underline{Hom}_\mathcal{A}(A_n,B_{n+r})$. We denote by $Hom(A,B)_r=\underset{n \in \mathbb{Z}}{\Pi}Hom_\A(A_n,B_{n+r})$ the underlying abelian group of maps of degree $r$. The tensor product is given : \begin{itemize}
    \item On objects by : \[(A\otimes B)_n:= \underset{k+l=n}{\bigoplus}A_k\otimes B_l\]
    \item On morphisms\footnote{Here, the morphims are the ones of the set-enriched version of $\mathrm{gr}\A$, i.e. the tensor product in $\mathrm{gr}\A$ is a priori not assumed to be $\A\mhyphen$enriched.}, the tensor product of $f$ and $g$ (both of degree $0$) is : \[(f\otimes g)_{|A_k\otimes B_l}:= f_k \otimes g_l\]
    \item The natural symmetry is given by : \[\tau^{A,B}_{|A_k\otimes B_l}:=(-1)^{kl}\tau^{A_k,B_l}_\mathcal{A}:A_k\otimes B_l \rightarrow B_l \otimes A_k  \]
\end{itemize} Here $\tau_\A$ denotes the symmetry of the tensor product of $\A$. \newline
The tensor product on graded objects and graded morphisms together with the internal object of morphisms $\underline{Hom}(A,B):= (\underline{Hom}(A,B)_r)_{r \in \mathbb{Z}}$ form a closed symmetric monoidal category with unit the unit object of $\A$ seen as a graded object concentrated in degree $0$. This category is still denoted $\mathrm{gr}\A$. This category, as a closed monoidal category, is enriched over itself. The category of graded objects of $\A$ is also $\A\mhyphen$enriched by only considering the degree $0$ part of the $\mathrm{gr}\A\mhyphen$enrichment of $\mathrm{gr}\A$, and the object of $\A$ of this enrichment is $\underline{Hom}^\A_{\mathrm{gr}\A}(X,Y):= \underset{n \in \mathbb{Z}}{\Pi}\underline{Hom}_\A(X_n,Y_n)$. 
\newline A chain complex $(V,d)$ in $\A$ is a $\mathbb{Z}\mhyphen$graded object $V$ of $\A$ together with a morphism, called a \textit{differential}, $d \in Hom_{\mathrm{gr}\A}(V,V)_{-1}$ of degree $-1$ which squares to zero : $d^2=0$. A morphism of chain complexes is a morphism of graded objects commuting with the differentials. This category, which is denoted $Ch(\A)$, is abelian but one can promote it to a closed symmetric monoidal category. The internal object of morphism is $\Homch((V,d_V),(W,d_W)):=(\Homgr(V,W),\partial)$. The differential $\partial: \Homgr(V,W)_{\bullet} \rightarrow \Homgr(V,W)_{\bullet-1}$ is given at each degree $r\in \mathbb{Z}$ by $\underline{Hom}(V,W)_r= \underset{n}{\Pi}\underline{Hom}_\mathcal{A}(V_n,W_{n+r})\rightarrow \underline{Hom}(V,W)_{r-1}=\underset{n}{\Pi}\underline{Hom}_\mathcal{A}(V_n,W_{n+r-1})$ by its value on the $i^{\text{th}}$ factor $\partial^i:\underset{n}{\Pi}\underline{Hom}_\mathcal{A}(V_n,W_{n+r}) \rightarrow \Hom(V_i,W_{i+r-1})$ as the adjoint map of the following composition : \[\begin{tikzcd}
{\underset{n}{\Pi}\underline{Hom}_\mathcal{A}(V_n,W_{n+r})} \arrow[r, "\partial^i"] \arrow[d, "p_{can}"', two heads]                     & {\underline{Hom}_\mathcal{A}(V_i,W_{i+r-1})} \\
{\underline{Hom}_\mathcal{A}(V_i,W_{i+r})\oplus\underline{Hom}_\mathcal{A}(V_{i-1},W_{i+r-1})} \arrow[ru, "\alpha \oplus\beta"'] &                                                       
\end{tikzcd}\]
Using the adjunction $ -\otimes V_i\dashv \Hom(V_i,-)$ and that tensor products and coproducts commute, ${\partial^i}^*$ is the composition : \[\begin{tikzcd}
({\underset{n}{\Pi}\underline{Hom}_\mathcal{A}(V_n,W_{n+r})}) \otimes V_i\arrow[r, "{\partial^i}^*"] \arrow[d, "p_{can}"', two heads]                     & W_{i+r-1} \\
({\underline{Hom}_\mathcal{A}(V_i,W_{i+r}) \otimes V_i)\oplus (\underline{Hom}_\mathcal{A}(V_{i-1},W_{i+r-1})}\otimes V_i) \arrow[ru, "(d_W\circ ev) \oplus (-1)^{r+1}(ev \circ (\mathrm{Id}\otimes d_V))"'] &                                                       
\end{tikzcd}\]
The tensor product of two chain complexes $(V,d_V)$ and $(W,d_W)$ has $V\otimes W$ as underlying graded object and the differential $d_{V\otimes W}$ is given on a summand $V_k\otimes W_l$ by :\[((d_V \otimes \mathrm{Id}) , (-1)^k\mathrm{Id}\otimes d_W): V_k\otimes W_l \rightarrow (V_{k-1}\otimes W_l) \oplus (V_k\otimes W_{l-1}) \] 
\begin{lemma}\label{enriched}The data of the above internal object of morphisms, tensor product and symmetry coming from the one on graded objects, define a closed symmetric monoidal structure on $Ch(\A)$. Moreover the canonical forgetful functor : \[U:Ch(\A) \rightarrow \mathrm{gr}\A\] is $\mathrm{gr}\A\mhyphen$enriched, thus $\A\mhyphen$enriched by pushing forward the enrichement along the (lax) monoidal functor $\mathrm{gr}\A \rightarrow \A$, $V_{\bullet} \mapsto V_0$.
\end{lemma}
In fact one can weaken this enrichement by only enriching $Ch(\A)$ over $\A$ by the same methods, we denote it $\underline{Hom}_{Ch(\A)}^\A(V,W)$, and this satisfies a natural adjunction formula for all $Y \in Ch(\A)$ : \[\begin{tikzcd}
\mathcal{A} \arrow[rr, "-\otimes Y", shift left=2,""{name=A, below}] &  & Ch(\mathcal{A}) \arrow[ll, "{\underline{Hom}^\mathcal{A}_{Ch(\mathcal{A})}(Y,-)}", shift left=2,""{name=B,above}]\ar[from=A, to=B, symbol=\dashv]
\end{tikzcd}\]  More explicitely $\underline{Hom}^\mathcal{A}_{Ch(\mathcal{A})}(Y,X)$ are the $0\mhyphen$cycles of the chain complex $\Homch(X,Y)$, that is \[\underline{Hom}^\mathcal{A}_{Ch(\mathcal{A})}(X,Y)\simeq Z_0(\Homch(X,Y)),\] representing the object of $\A$ of chain complexes morphisms from $X$ to $Y$.
\begin{example}
    If $\A=\mathrm{Mod}_k$ is the category of $k\mhyphen$modules over the commutative ring $k$, then for all $X,Y \in Ch(\A)$ the object $\underline{Hom}^\mathcal{A}_{Ch(\mathcal{A})}(X,Y)$ is the $k\mhyphen$module of chain complex morphisms from $X$ to $Y$.
\end{example}
\begin{lemma}
    For each $n\in \mathbb{Z}$, the $n^{\text{th}}\mhyphen$projection map :\[\begin{array}{ccccc}
ev_n & : & Ch(\A) & \to & \A \\
 & & V & \mapsto & V_n \\
\end{array}\] defines an $\A\mhyphen$enriched functor.
\end{lemma}
\begin{demo}
    Let us build natural maps $\underline{Hom}^\A_{Ch(\A)}(V,W) \rightarrow \underline{Hom}_\A(V_n,W_n)$ for every $V,W \in Ch(\A)$. This map is the adjoint map of the following composition (the dotted arrow): \[\begin{tikzcd}
{\underline{Hom}^\mathcal{A}_{Ch(\mathcal{A})}(V,W) \otimes V_n} \arrow[rr, dotted] \arrow[d, "U\otimes \mathrm{Id}"'] &                                                               & W_n \\
{\underline{Hom}^\A_{gr\mathcal{A}}(V,W)\otimes V_n} \arrow[d, "\simeq"']                                               &                                                               &     \\
{\underset{j}{\Pi}\underline{Hom}_\mathcal{A}(V_j,W_j)\otimes V_n} \arrow[r, "p_{can}"', two heads]                    & {\underline{Hom}_\mathcal{A}(V_n,W_n)\otimes V_n} \arrow[ruu, "ev"'] &    
\end{tikzcd}\] Then it is routine to check that this defines an $\A\mhyphen$enriched functor.
\end{demo}

This allows us to define a new endormorphism operad for chain complexes :
\begin{definition}
    The $\A\mhyphen$enriched endomorphism operad of a chain complex $(V,d)$, denoted $End^\A(V,d)$ or simply $End^\A(V)$ if no confusion is possible, is the $\A\mhyphen$operad with $n\mhyphen$ary operations $End^\A(V)_n:=\underline{Hom}_{Ch(\A)}^\A((V,d)^{\otimes n},(V,d))$. The operadic structure is given by the family of maps, for all $l\in \mathbb{N}$ and every $n_1,...,n_l \in \mathbb{N}$ : $End^\A(V)_l\otimes End^\A(V)_{n_1}\otimes ... \otimes End^\A(V)_{n_l}$ as the adjoint map of the dotted arrow  : \[\begin{tikzcd}
{\underline{Hom}_{Ch(\mathcal{A})}^\mathcal{A}(V^{\otimes l},V)\otimes (\underset{1 \leq i \leq l}{\otimes }\underline{Hom}_{Ch(\mathcal{A})}^\mathcal{A}(V^{\otimes n_i},V))\otimes V^{\otimes \underset{i}{\Sigma}n_i}} \arrow[d, "\mathrm{Id} \otimes \underline{\otimes}\otimes \mathrm{Id}"'] \arrow[r, dotted] & V                                                                                                      \\
{\underline{Hom}_{Ch(\mathcal{A})}^\mathcal{A}(V^{\otimes l},V)\otimes \underline{Hom}_{Ch(\mathcal{A})}^\mathcal{A}(V^{\otimes \underset{i}{\Sigma}n_i},V^{\otimes l})\otimes V^{\otimes \underset{i}{\Sigma}n_i}} \arrow[r, "\mathrm{Id} \otimes ev"]                 & {\underline{Hom}_{Ch(\mathcal{A})}^\mathcal{A}(V^{\otimes l},V)\otimes V^{\otimes l}} \arrow[u, "ev"']
\end{tikzcd}\] Here we denote $\underline{\otimes}$ the morphism induced by \[\underline{Hom}_{Ch(\mathcal{A})}^\mathcal{A}(X,Y) \otimes \underline{Hom}_{Ch(\mathcal{A})}^\mathcal{A}(Z,T) \rightarrow \underline{Hom}_{Ch(\mathcal{A})}^\mathcal{A}(X\otimes Z,Y\otimes T)\] for all $X,Y,Z,T \in Ch(\A)$.
\end{definition}
The canonical monoidal inclusion $\A\rightarrow Ch(\A)$ sending an object of $\A$ to the chain complex concentrated in degree $0$, induces a fully faithful inclusion $\iota: Op(\A) \rightarrow Op(Ch(\A))$ which gives us the following :
\begin{lemma}\label{repalg}
    For every chain complex $(V,d)$ and every operad $\p$ of $\A$, there is a natural isomorphism : \[\{\iota(\p)\mhyphen\text{alg structure on } (V,d)\} \simeq Hom_{Op(\A)}(\p,End^\A(V,d))\]
\end{lemma}
\begin{remarque}\label{nfoldchain}
As a consequence, there is a closed monoidal structure on the category of chain complexes of a closed monoidal category\footnote{In the rest of this paper we will call closed monoidal any closed symmetric monoidal abelian category. Hence, any closed monoidal category $\A$ gives rise to closed monoidal categories $Ch(\A)$ and $nCh(\A)$.}. In particular the category of $n\mhyphen$fold chain complexes $nCh(\A)$ of $\A$ is again closed monoidal. By a direct induction of the previous process we get an $\A\mhyphen$enrichement of the category $nCh(\A)$ which allows us to represent $\iota(\p)\mhyphen$algebra structures on any $n\mhyphen$fold chain complex $V \in nCh(\A)$ and for all $\A\mhyphen$operad $\p$ : \[\{\iota(\p)\mhyphen\text{alg structure on }V\}\simeq Hom_{Op(\A)}(\p, End^\A(V))\]
\end{remarque}
From now on, we drop off the "$\iota$" in the notation $\iota(\p)$, and $\p$ is as always a $\A\mhyphen$operad.
%\begin{rappel}
    %The category $Ch(\A)$ of chain complexes over $\A$ is abelian monoidal closed and the canonical inclusion $V \mapsto (...0\rightarrow V\rightarrow 0 ...) $, which sends an object of $\A$ to the chain complexe concentrated in degree $0$, is monoidal. The tensor product of two chain complexes $C_\bullet$ and $D_\bullet$ is given by :  \[(C_\bullet \otimes D_\bullet )_n :=\underset{k+l=n}{\bigoplus}C_k \otimes D_l\] with differentials on each summand : \[(d\otimes Id) \oplus ((-1)^{k}Id\otimes d):C_k\otimes D_l \rightarrow (C_{k-1}\otimes D_l) \oplus (C_k \otimes D_{l-1}) \]The symmetry $\tau_{Ch(\mathcal{A})}$ is given by the Koszul sign rule : if $C,D \in Ch(\mathcal{A})$ $\tau_{Ch(\mathcal{A})} : C \otimes D \rightarrow D \otimes C $ is induced by its value on each summand : $(-1)^{kl}\tau_{C,D}^\mathcal{A}:C_k\otimes D_l \rightarrow D_l \otimes C_k $. The tensor product of two chain map $f:C\rightarrow D$ and $g:C'\rightarrow D'$ is given by : 
%\end{rappel}

\subsection{Universal semi-direct product, derived operations and codescent}
Throughout this section $V$ is a chain complex in $\A$ concentrated in degrees $0$ and $1$, that is $V=...0\rightarrow V_1 \overset{d}{\rightarrow} V_0\rightarrow 0 ...$. The goal here is to construct a morphism of symmetric $\A\mhyphen$operads :\[\rtimes : End^\A(V) \rightarrow End(V_1\oplus V_0)\] which gives a universal semi-direct product in the sense that, for any operad $\p$ and any $\p\mhyphen$algebra structure on $V$, say $\gamma : \p \rightarrow End^\A(V)$, $\gamma$ gives rise by post-composing it with $\rtimes$ to a $\mathrm{Cat}^1\mhyphen\p\mhyphen$algebra, \[\begin{tikzcd}
V_1\rtimes V_0 \arrow[r, "0\oplus \mathrm{Id}", shift left=2] \arrow[r, "d \oplus \mathrm{Id}"', shift right=2] & V_0 \arrow[l, "i_2" description]
\end{tikzcd}\]We first recall some standard results to construct a map $\partial : (V_1\oplus V_0)^{\otimes n}\rightarrow (V^{\otimes n})_1 \oplus (V^{\otimes n})_0$.
\begin{definition}
Let $p,q$ be integers. A $(p,q)\mhyphen$\textit{shuffle} is a permutation $\sigma \in \Sigma_{p+q}$ such that $
\sigma(1) < ... < \sigma(p) $ and $
\sigma(p+1)< ... < \sigma(q)$. We denote by Sh$(p,q)$ the set of $(p,q)\mhyphen$shuffles.
\end{definition}
\begin{remarque}
By abuse of notation a $(p,q)\mhyphen$shuffle is represented as a pair $(\mu,\nu)$, where $\{\mu,\nu\}$ is a partition of the set $\{1,2,...,p+q\}$ with $\mu=\{\mu_1,...,\mu_p\}$ has $p$ elements and $\nu=\{\nu_1,...,\nu_q\}$ has $q$ elements.
\end{remarque}
Given $(\mu,\nu)$ a $(p,q)\mhyphen$shuffle, we denote $V^{(\mu,\nu)}$ the tensor product of copies of $V_1$ and $V_0$ where $V_1$ appears in places $\mu_1,...,\mu_p$ and $V_0$ appears in places $\nu_1,...,\nu_q$.
\begin{example}
If $(\mu,\nu)=(\{1,3\}, \{2\})$, then $V^{(\mu,\nu)}=V_1\otimes V_0 \otimes V_1$.
\end{example}
\begin{lemma}
  There is the following decomposition :\[(V_1\oplus V_0)^{\otimes n}=\underset{p+q=n}{\bigoplus} \underset{(\mu,\nu) \in \text{Sh}(p,q)}{\bigoplus}V^{(\mu,\nu)}\]  
\end{lemma}
If $p \geq 2$ and $\mu \neq \emptyset$, then for all $(\mu,\nu) \in \mathrm{Sh}(p,q)$, the differential $d:V_1 \rightarrow V_0$ induces a morphism $\partial: V^{(\mu,\nu)}\rightarrow V^{(\mu_1, \{1,...,p+q\}\setminus {\mu_1})} $ by applying $d$ on each copy of $V_1$ except on the left most one\footnote{We could have chosen to apply $d$ on every copies of $V_1$ except on exactly one, this does not matter for Definition \ref{defproduitsemi} because of the action of the symmetric groups.}, and applying the identity on copies of $V_0$. We choose to extend $\partial$ by the identity on $(0,n)\mhyphen$shuffles and on $(1,n-1)\mhyphen$shuffles. Thus we get a morphism :\[\partial : (V_1 \oplus V_0)^{\otimes n} \rightarrow (V^{\otimes n})_1\oplus (V^{\otimes n})_0\]
\begin{definition}\label{defproduitsemi}
    The morphism :\[\rtimes : End^\A(V) \rightarrow End(V_1 \oplus V_0) \]is given, for each $n$, by the dotted arrow on the $n\mhyphen$ary operations, which is the composition of the other ones : \[\begin{tikzcd}
{\underline{Hom}^\A_{Ch(\A)}(V^{\otimes n}, V)} \arrow[d, "{(ev_1,ev_0)}"'] \arrow[r, "\rtimes", dotted]                & {\underline{Hom}_\A((V_1\oplus V_0)^{\otimes n}, V_1 \oplus V_0)}                                        \\
{\underline{Hom}_\A((V^{\otimes n})_1, V_1) \oplus \underline{Hom}_\A((V^{\otimes n})_0, V_0) } \arrow[r, "diag"'] & {\underline{Hom}_\A((V^{\otimes n})_1\oplus (V^{\otimes n})_0, V_1 \oplus V_0)} \arrow[u, "\partial^*"']
\end{tikzcd}\] 
\end{definition}
\begin{proposition}\label{codes}
    $\rtimes : End^\A(V) \rightarrow End(V_1 \oplus V_0)$ is a morphism of $\A\mhyphen$operads.
\end{proposition}
\begin{demo}
    %For the sake of simplicity, let us assume that $\A$ is the category of $k\mhyphen$modules for some ring $k$. We remark for the careful reader that it is not a use of the Freyd-Mitchell embedding theorem because this theorem does not necessarily provide a commutative ring, in particular there is no reason for our tensor product in $\A$ to be the tensor product over a commutative ring $\underset{k}{\otimes}$. However the only properties of the category of $k\mhyphen$modules we use are the fact that the tensor product commutes with colimits and that the set of maps are abelian groups. These properties are true in $\A$, in particular the careful reader can transpose this set-based proof into a diagrammatic proof to get a complete proof for any $\A$. \newline 
    First we prove that $\rtimes$ is a morphism of $\Sigma_*\mhyphen$modules\footnote{A $\Sigma_*\mhyphen$module in a category $\mathcal{C}$ is a functor $\Sigma=\underset{n \in \mathbb{N}}{\amalg}\Sigma_n \rightarrow \mathcal{C}$. The data of such a functor is equivalent to the data of a collection of objects $M(n)\in \mathcal{C}$ for all $n \in \mathbb{N}$, where each $M(n)$ is endowed with a right action of the symmetric group $\Sigma_n$. }\cite{fresse2007modules}, that is for all $n \in \mathbb{N},\sigma \in \Sigma_n$ we have $\rtimes \circ \sigma = \sigma \circ \rtimes : \underline{Hom}^\A(V^{\otimes n}, V) \rightarrow \underline{Hom}((V_1 \oplus V_0)^{\otimes n}, V_1\oplus V_0)$. In fact we prove that the difference of the adjoint maps is zero :\[\rtimes \circ \sigma^* -\sigma \circ \rtimes^* =0: \underline{Hom}^\A(V^{\otimes n}, V) \otimes (V_1 \oplus V_0)^{\otimes n} \rightarrow V_1 \oplus V_0 \]The left side summand $\rtimes \circ \sigma^*$ is given by the action of $\sigma $ on $V^{\otimes n}$ and applying $\partial $ on $(V_1\oplus V_0)^{\otimes n}$, then applying the evaluation maps. The right hand side is given by acting with $\sigma $ on $(V_1 \oplus V_0)^{\otimes n}$, applying $\partial$ on it, then using evaluation maps. That is, one has to prove that the following diagram commutes : \begin{equation}\begin{tikzcd}
{\underline{Hom}^\A(V^{\otimes n},V)\otimes (V_1 \oplus V_0)^{\otimes n}} \arrow[r] \arrow[rd, "0"'] & {\underline{Hom}^\A(V^{\otimes n},V)\otimes ((V^{\otimes n})_1 \oplus (V^{\otimes n})_0)} \arrow[d, "{(ev_1\circ(\mathrm{Id}-\sigma),ev_0\circ(\mathrm{Id}-\sigma))}"] \\
                                                                                                                                                            & V_1\oplus V_0                                                                                                                                       
\end{tikzcd}\label{eq:demop} \tag{2}\end{equation} Here the upper arrow is given by $\mathrm{Id}\otimes(\partial \circ \sigma)-\mathrm{Id}\otimes \partial$. As $V$ is concentrated in degrees $0$ and $1$, the following diagram commutes : \[\begin{tikzcd}
{\underline{Hom}^\A(V^{\otimes n},V)} \arrow[r, "ev_1"] \arrow[rd, "0"'] & {\underline{Hom}^\A((V^{\otimes n})_1,V_1)} \arrow[d, "d_{V^{\otimes n}}^*"] \\
                                                                         & {\underline{Hom}^\A((V^{\otimes n})_2,V_1)}                                 
\end{tikzcd}\]So, as the composition of the two maps of the diagram \eqref{eq:demop} factors through the composition of the previous diagram, we are done. Now we have to prove that $\rtimes$ is compatible with the operadic structures, that is it commutes with partial compositions $\circ_i$. Thus we have to prove that for all $n,m \in \mathbb{N}$ and for all $i \in \llbracket 1,n\rrbracket$, the following diagram is commutative : \[\begin{tikzcd}
    End^\A(V)_n\otimes End^\A(V)_m \arrow[d, "\circ_i"'] \arrow[r, "\rtimes \otimes \rtimes"] & End^\A(V_1 \oplus V_0)_n\otimes End^\A(V_1 \oplus V_0)_m \arrow[d, "\circ_i"] \\
    End^\A( V)_{n+m-1} \arrow[r, "\rtimes"']                                                  & End^\A(V_1 \oplus V_0)_{n+m-1}                                               
     \end{tikzcd} \]As partial compositions in endomorphism operads are given by precomposition and as $\rtimes$ is also given by precomposition with $\partial$, then the diagram is commutative because $d:V_1 \rightarrow V_0$ is a derivation for the $End^\A(V)\mhyphen$algebra structure on $V$. 
\end{demo}
\begin{proposition}\label{semitocat}
    For any $\p\mhyphen$algebra structure on $V$, the $\p\mhyphen$algebra structure on $V_1\oplus V_0$ given by the operad morphism $\p\rightarrow End^\A(V)\overset{\rtimes}{\rightarrow}End(V_1 \oplus V_0)$ induces a $\mathrm{Cat}^1\mhyphen\p\mhyphen$algebra structure on the reflexive graph $\begin{tikzcd}
V_1\oplus V_0 \arrow[r, "0\oplus \mathrm{Id}", shift left=2] \arrow[r, "d \oplus \mathrm{Id}"', shift right=2] & V_0 \arrow[l, "i_2" description]
\end{tikzcd}$.
\end{proposition}
\begin{demo}
We have to prove that it is a reflexive graph of $\p\mhyphen$algebras and that the $\mathrm{Cat}^1\mhyphen$condition \ref{cat^1} holds. It is clear that $V_0 \overset{i}{\hookrightarrow} V_1 \oplus V_0$ and $V_1 \oplus V_0 \overset{0\oplus \mathrm{Id}}{\longrightarrow} V_0$ are $\p\mhyphen$algebra morphisms. In particular $V_1=ker(0\oplus \mathrm{Id}:V_1\oplus V_0 \rightarrow V_0)$ is naturally endowed with a $\p\mhyphen$algebra structure as a kernel of a $\p\mhyphen$algebra morphism. Let us check that $d\oplus \mathrm{Id} $ is a $\p\mhyphen$algebra morphism, that is the following diagram is commutative : \begin{tikzcd}
\p_n\otimes (V_1 \oplus V_0)^{\otimes n} \arrow[d, "\mu^{V_1\oplus V_0}"'] \arrow[r, "\mathrm{Id}\otimes (d\oplus \mathrm{Id})^{\otimes n}"'] & \p_n\otimes V_0^{\otimes n} \arrow[d, "\mu^{V_0}"] \\
V_1\oplus V_0 \arrow[r, "d \oplus \mathrm{Id}"']                                                                                    & V_0                                               
\end{tikzcd} where $\mu^{V_1\oplus V_0}$ is the structure map induced by the adjoint map of $\p_n\rightarrow End(V_1\oplus V_0)_n$. Recall that $(V_1\oplus V_0)^{\otimes n}=\underset{p+q=n}{\bigoplus} \underset{(\mu,\nu) \in Sh(p,q)}{\bigoplus}V^{(\mu,\nu)}$, so as $\otimes$ commutes with $\oplus$, it is enough to ckeck the commutativity of the previous diagram on summands $\p_n\otimes V^{(\mu,\nu)}$. If $V^{(\mu,\nu)}=V_0\otimes...\otimes V_0$, then this is obvious. If $V^{(\mu,\nu)}$ is not of this form, then $\p_n\otimes V^{(\mu,\nu)} \rightarrow V_1\oplus V_0$ restricts to $\p_n\otimes V^{(\mu,\nu)} \overset{\mu^{V^1}}{\rightarrow} V_1$ by definition of the structure map of $V_1\oplus V_0$. So $\mu^{V_1}$ fits into a commutative diagram  \[\begin{tikzcd}
{\p_n\otimes V^{(\mu,\nu)}} \arrow[r, "\mu^{V_1}"] \arrow[d, "\mathrm{Id}\otimes \partial"'] & V_1 \\
\p_n\otimes (V^{\otimes n})_1 \arrow[ru, "\mu^V"']                                  &    
\end{tikzcd}\] As the following diagram commutes : \begin{tikzcd}
\p_n\otimes (V^{\otimes n})_1 \arrow[d, "\mu^V"'] \arrow[r, "\mathrm{Id}\otimes d"] & \p_n\otimes V^{\otimes n}_0 \arrow[d, "\mu^V=\mu^{V_0}"] \\
V_1 \arrow[r, "d"']                                                        & V_0                                                     
\end{tikzcd} we are left to prove that \begin{tikzcd}
{\p_n\otimes V^{(\mu,\nu)}} \arrow[d, "\mathrm{Id}\otimes \partial"'] \arrow[r, "\mathrm{Id} \otimes \bar{d}"] & \p_n\otimes V^{\otimes n}_0 \\
\p_n\otimes (V^{\otimes n})_1 \arrow[ru, "\mathrm{Id}\otimes d_{V^{\otimes n}}"']                     &                            
\end{tikzcd} commutes, where $\bar{d}: V^{(\mu,\nu)}\rightarrow V_0^{\otimes n}$ is the map which applies $d$ on each copy of $V_1$. But this last diagram is clearly commutative as $d_{V^{\otimes n}}\circ \partial=\bar{d}$. To prove that the $\mathrm{Cat}^1\mhyphen$condition is fulfilled, according to Proposition \ref{peiffer=Cat^1}, it is enough to prove that $V_1 \overset{d}{\rightarrow}V_0$ is a crossed module, which is the purpose of the following lemma.
\end{demo}
\begin{lemma}
    For the induced structure of $\p\mhyphen$algebra on $V_1$ (as the kernel of the $\p\mhyphen$algebra morphism $0\oplus \mathrm{Id}:V_1\oplus V_0 \rightarrow V_0$), and for the canonical internal action of $V_0$ on $V_1$, the precrossed module $V_1 \overset{d}{\rightarrow}V_0$ satisfies the Peiffer identity, that is, it is a crossed module of $\p\mhyphen$algebras.
\end{lemma}
\begin{demo}
    Here, the Peiffer identity (\ref{eq:peiffer}) is expressed by the commutativity of the following diagram :\[\begin{tikzcd}
\p_n\underset{}{\otimes} V^{\otimes i}_1\otimes V^{\otimes n-i}_0 \arrow[d, "\mathrm{Id}\otimes \mathrm{Id} \otimes d^{i-1}\otimes \mathrm{Id}"'] \arrow[r] & \p\circ (V_0| V_1) \arrow[rd,"\rho"] \\
\p_n\underset{}{\otimes} V_1\otimes V^{\otimes n-1}_0 \arrow[r] & \p\circ (V_0| V_1) \arrow[r, "\rho"]& V_1                                                        &                  
\end{tikzcd}\] It is ensured by the map $\rtimes : End^\A(V) \rightarrow End(V_1 \oplus V_0)$, as it is given by precomposition with $\partial$. Indeed, on the summand $V^{\otimes i}_1\otimes V^{\otimes n-i}_0$, the map $\partial$ is given by $\mathrm{Id}\otimes d^{\otimes i-1}\otimes \mathrm{Id}$.\end{demo}

\vspace{1cm}
We prove the converse now.
\begin{proposition}
    Given a crossed module of $\p\mhyphen$algebras $(V_0 \overset{\rho}{\curvearrowright}V_1, d : V_1 \rightarrow V_0)$, the family of maps :\begin{equation} \mu^k_n:(\p_n\underset{\Sigma_n}{\otimes}V^{\otimes n})_k\rightarrow V_k , \hspace{1cm}\forall n,k \in \mathbb{N} \end{equation} define a $\p\mhyphen$algebra structure on the complex $V=...0\rightarrow  V_1 \overset{d}{\rightarrow} V_0$, where : \begin{itemize}
        \item $\mu^0_n: (\p_n\underset{\Sigma_n}{\otimes}V^{\otimes n})_0 = \p_n\underset{\Sigma_n}{\otimes}V_0^{\otimes n} \rightarrow V_0$ for $n \geq 0$, $\mu^0_n$ is the structure map of the $\p\mhyphen$algebra structure of $V_0$, that is $\mu^0=\mu^{V_0}$.
        \item  $\mu^1_n: (\p_n\underset{\Sigma_n}{\otimes}V^{\otimes n})_1 \rightarrow V_1$ are the structure maps given by the action of $V_0$ on $V_1$. More precisely it is induced by the structure of $V_0\mhyphen\p\mhyphen$module on $V_1$.        
       \item $\mu^{k \geq 2}_n=0$
    \end{itemize}This construction provides a functor $\mathrm{XMod}(\p\mhyphen alg)\rightarrow \p\mhyphen alg_{\leq 1}(Ch(\A))$ from the category of crossed modules of $\p\mhyphen$algebras (Definition \ref{Peiffer}) to the category of differential graded $\p\mhyphen$algebras concentrated in degrees $0$ and $1$.
\end{proposition}
\begin{demo}
    We first prove that the family of maps $(\mu_n^k)_{k \in \mathbb{N}}$ defines chain complexes morphisms, that is for all $n$, the following diagram commutes :\[\begin{tikzcd}
... \arrow[r] \arrow[d] & \p_n\underset{\Sigma_n}{\otimes}(V^{\otimes n})_2 \arrow[d, "0"] \arrow[r, "\mathrm{Id}\otimes d_{V^{\otimes n}}"] & \p_n\underset{\Sigma_n}{\otimes}(V^{\otimes n})_1 \arrow[d, "\mu^1_n"] \arrow[r, "\mathrm{Id}\otimes d_{V^{\otimes n}}"] & \p_n\underset{\Sigma_n}{\otimes}V_0^{\otimes n} \arrow[d, "\mu^0_n"] \\
... \arrow[r]           & 0 \arrow[r, "0"']                                                                                         & V_1 \arrow[r, "d"']                                                                                             & V_0                                                                 
\end{tikzcd}\]The right hand square commutes because of equivariance of $d$ with respect to the action of $V_0$ on itself and on $V_1$. The middle square commutes because of the Peiffer identity \ref{Peiffer} for $i=2$, $X=V_1$ and $B=V_0$. The left hand squares trivially commutes as $0$ is a terminal object. Moreover, these morphisms are, by definition, compatible with the actions of symmetric groups, that is, for all $n \in \mathbb{N}$, the adjoint maps $\p_n\rightarrow End^\A(V)_n$ are equivariant with respect to the action of the symmetric group $\Sigma_n$. It remains to prove that the family $\mu_\bullet$ defines a $\p\mhyphen$algebra structure on $V$, thus we are led to prove that it is unital and associative. Unitality is a direct consequence of the definition. To prove it is associative, we have to prove for all $n,m \in \mathbb{N}$ and for all $i \in \llbracket 1,n\rrbracket$, that the following diagram is commutative :\[\begin{tikzcd}
    \p_n\otimes\p_m\otimes V^{\otimes n+m-1} \arrow[d, "\circ_i\otimes \mathrm{Id}"'] \arrow[r, "\mathrm{Id}\otimes \mu_m\otimes\mathrm{Id}^{\otimes n-i}_{V}"'] & \p_n\otimes V^{\otimes n} \arrow[d, "\mu_n"] \\
    \p_{n+m-1}\otimes V^{\otimes n+m-1} \arrow[r, "\mu_{n+m-1}"']                                                                                               & V                                           
    \end{tikzcd}\] 

Here we denote $\circ_i$ the partial composition induced by the structure map of the operad $\p$. As $V$ is a chain complex concentrated in degrees $0$ and $1$, the commutativity of this diagram is equivalent to the commutativity of its degree $0$ part and its degree $1$ part. For the degree $0$ part, it is the one induced by restricting the top left corner to $ \p_n\otimes\p_m\otimes V_0^{\otimes n+m-1}$, which is commutative as the family $\mu^0_\bullet$ is, by definition, the structure map of the $\p\mhyphen$algebra $V_0$, thus associativity is fulfilled. For the degree $1$ part, it is commutative because the family $\mu^1_\bullet$ is, by definition, given by structure map  of the $V_0\mhyphen\p\mhyphen$module structure on $V_1$, and the associativity of the action\footnote{Here, we mean the $V_0\mhyphen\p\mhyphen$module structure on $V_1$, that is the action of $V_0$ on the abelian $\p\mhyphen$algebra $V_1$.} is exactly expressed by the commutativity of this diagram.

\end{demo}\newline

As a corollary, we get the following theorem : 
\begin{theorem}\label{equivcrossedmod}
    The previous constructions define equivalences of categories : \[\begin{tikzcd}
                                              & \mathrm{Cat}(\p\mhyphen alg) \arrow[ld] \arrow[rd] &                                                             \\
\mathrm{XMod}(\p\mhyphen alg) \arrow[ru, shift left=2] &                                                   & \p\mhyphen alg_{\leq 1}(Ch(\A)) \arrow[lu, shift right=2]
\end{tikzcd}\]Here the map $\mathrm{Cat}(\p\mhyphen alg) \rightarrow \p\mhyphen alg_{\leq 1}(Ch(\A))$ is the composition :\[\mathrm{Cat}(\p\mhyphen alg) \rightarrow \mathrm{XMod}(\p\mhyphen alg) \rightarrow \p\mhyphen alg_{\leq 1}(Ch(\A))\] induced by the previous proposition.
\end{theorem}
The composition $\p\mhyphen alg_{\leq 1}(Ch(\A))\rightarrow \mathrm{Cat}(\p\mhyphen alg) \rightarrow \mathrm{XMod}(\p\mhyphen alg)$ is called the derived operations functor in reference to \cite{Kosmann_Schwarzbach_2004}, or the codescent structure functor in the sense that for a global $\p\mhyphen$algebra structure on the object $V=...0\rightarrow V_1 \rightarrow V_0$ we get local $\p\mhyphen$algebra structures on $V_1$ and $V_0$. The other composition :\[\mathrm{XMod}(\p\mhyphen alg) \rightarrow \mathrm{Cat}(\p\mhyphen alg) \rightarrow \p\mhyphen alg_{\leq 1}(Ch(\A))\] is dually a descent structure functor. So the "local" notion of crossed module, namely the category $\mathrm{XMod}(\p\mhyphen alg)$ is equivalent to the "global" one, namely $\p\mhyphen alg_{\leq 1}(Ch(\A))$.\newline
\begin{demo}
    We are left to prove that the pair of functors \begin{tikzcd}
\mathrm{XMod}(\p\mhyphen alg) \arrow[r, shift left] &  \p\mhyphen alg_{\leq 1}(Ch(\A))\arrow[l, shift left]
\end{tikzcd} are inverse equivalences. On the one hand, the composition \[\mathrm{XMod}(\p\mhyphen alg) \rightarrow \p\mhyphen alg_{\leq 1}(Ch(\A)) \rightarrow \mathrm{XMod}(\p\mhyphen alg)\] is isomorphic to the identity because of the Peiffer condition which fully determines the $\p\mhyphen$algebra structure on $V_1$ as long as we have the $V_0\mhyphen\p\mhyphen$module structure and the $V_0$ equivariant map $d$. The other composition is clearly isomorphic to the identity.
\end{demo}
\begin{example}
    Assume $\p=\mathfrak{L}ie$ is the operad encoding Lie algebras and $\A=Vect_\mathbb{K}$ is the category of $\mathbb{K}\mhyphen$vector spaces. Then if $g=...0\rightarrow g_1\rightarrow g_0$ is a differential graded Lie algebra, the derived bracket on $g_1$ is given by $[x,y]_{g_1}:=[x,dy]_g$, and this leads to a Lie algebra structure on $g_1$. More generally, if $\mu \in \p_n$ and if $V=...0\rightarrow V_1\rightarrow V_0$ is a differential graded $\p\mhyphen$algebra, then the derived structure on $V_1$ is given by $\mu^{V_1}(x_1,...,x_n):=\mu^V(x_1,dx_2,...,dx_n)$.
\end{example}
\subsection{ $n\mhyphen$crossed$\mhyphen\p\mhyphen$algebras}
The category of $n\mhyphen$crossed modules of $\p\mhyphen$algebras can be recursively defined, using the local notion of crossed modules, as $\mathrm{X}^{n+1}\mathrm{Mod}(\p\mhyphen alg)=\mathrm{XMod(X}^n\mathrm{Mod}(\p\mhyphen alg))$ and it is equivalent to $\mathrm{Cat}^n(\p\mhyphen alg)$ the category of $n\mhyphen$fold categories internal to $\p\mhyphen$algebras. The problem of this definition is that such objects are quite involved. The notion of $n\mhyphen$crossed modules is better understood using the global definition of crossed modules. 
\begin{definition}
    An $n\mhyphen$crossed$\mhyphen\p\mhyphen$algebra is a $\p\mhyphen$algebra structure on an $n\mhyphen$fold chain complex $C_{\bullet,...,\bullet}$ concentrated in degrees $\epsilon_1,...,\epsilon_n$ where $\epsilon_i \in \{0,1\}$. A morphism of $n\mhyphen crossed\mhyphen\p\mhyphen$algebras is a morphism of $\p\mhyphen$algebras. They form a category denoted $\p\mhyphen alg_{\leq 1}(nCh(\A))$.
\end{definition}
\begin{example}\begin{itemize}
    \item A $2\mhyphen$crossed$\mhyphen\p\mhyphen$algebra is a $\p\mhyphen$algebra structure on a square :\[ \begin{tikzcd}
0 \arrow[d, dotted] & 0 \arrow[d, dotted]           &                     \\
{C_{0,1}} \arrow[d] & {C_{1,1}} \arrow[d] \arrow[l] & 0 \arrow[l, dotted] \\
{C_{0,0}}           & {C_{1,0}} \arrow[l]           & 0 \arrow[l, dotted]
\end{tikzcd} \]
\item Given $g=(g_0 \leftarrow g_1 \leftarrow 0 ...)$ a differential graded Lie algebra, denoting $\mathrm{Der}(g)_n$ the vector space of degree $n$ derivations of $g$, then the commutative square : \[\begin{tikzcd}
\mathrm{Der}(g)_1 \arrow[d] & g_1 \arrow[d] \arrow[l, "ad_\bullet"] \\
  \mathrm{Der}(g)_0        &  g_0 \arrow[l, "ad_\bullet"]          
\end{tikzcd}\] is naturally endowed with a $2\mhyphen crossed\mhyphen Lie\mhyphen$algebra structure. More generally, if $g$ is an $n\mhyphen crossed \mhyphen Lie\mhyphen$algebra, then $g \overset{ad}{\rightarrow} \mathrm{Der}(g)$ naturally inherits a structure of an $(n+1)\mhyphen crossed \mhyphen Lie\mhyphen$algebra. 

\end{itemize}
\end{example}
\begin{theorem}
    The categories $\p\mhyphen alg_{\leq 1}(nCh(\A))$ and $\mathrm{Cat}^n(\p\mhyphen alg)$ are equivalent.
\end{theorem}
\begin{demo} We prove this by induction on $n$.
    For $n=1$ it is the purpose of Theorem \ref{equivcrossedmod}. Now assume it is true for some $n$ and let us prove that $\p\mhyphen alg_{\leq 1}((n+1)Ch(\A))\simeq \mathrm{Cat}^{n+1}(\p\mhyphen alg)$. Any $C=C_{\bullet,...,\bullet} \in \p\mhyphen alg_{\leq 1}((n+1)Ch(\A))$ lies in $Ch(\A')$, where $\A':=nCh(\A)$, and $C=...0\rightarrow C_{1,\bullet,...,\bullet}\overset{d}{\rightarrow} C_{0,\bullet,...,\bullet}$ is viewed as a chain complex concentrated in degrees $0$ and $1$ in the category $\A'$. The $\p\mhyphen$algebra structure on $C$ is exactly a $dg\mhyphen\p'\mhyphen alg_{\leq 1}$ structure on $C$, where $\p'$ is the operad $\p$ considered as an operad of $\A'$. According to  Theorem \ref{equivcrossedmod}, $C$ gives rise to a category internal to the category of $\p'\mhyphen$algebras :\[\begin{tikzcd}
{C_{1,\bullet,...,\bullet}\rtimes C_{0,\bullet,...,\bullet}} \arrow[r, "0\oplus \mathrm{Id}", shift left=2] \arrow[r, "d\oplus \mathrm{Id}"', shift right=2] & {C_{0,\bullet,...,\bullet}} \arrow[l]
\end{tikzcd}\] This construction, as well as for the morphisms, produces an equivalence of categories between $\p\mhyphen alg_{\leq 1}((n+1)Ch(\A))$ and $\mathrm{Cat}(\p\mhyphen alg_{\leq 1}(Ch(\A)))$. Moreover $C_{1,\bullet,...,\bullet}\rtimes C_{0,\bullet,...,\bullet}$ and $C_{0,\bullet,...,\bullet}$ are themselves $n\mhyphen crossed \mhyphen\p\mhyphen$algebras, thus we get an equivalence of categories $\p\mhyphen alg_{\leq 1}((n+1)Ch(\A))\overset{\simeq}{\rightarrow}\mathrm{Cat}(\p\mhyphen alg_{\leq 1}(nCh(\A)))$. By induction hypothesis, the latter is equivalent to $\mathrm{Cat}(\mathrm{Cat}^n(\p\mhyphen alg))$. Thus, there is a sequence of equivalences of categories  \[\p\mhyphen alg_{\leq 1}((n+1)Ch(\A))\overset{\simeq}{\rightarrow}\mathrm{Cat}(\p\mhyphen alg_{\leq 1}(nCh(\A)))\overset{\simeq}{\rightarrow} \mathrm{Cat}^{n+1}(\p\mhyphen alg) \]
\end{demo}
\begin{example}
    If $g = \begin{tikzcd}
{g_{0,1}} \arrow[d, "d^v_0"] & {g_{1,1}} \arrow[d,"d^v_1"] \arrow[l,"d^h_1"] \\
{g_{0,0}}           & {g_{1,0}} \arrow[l,"d^h_0"]          
\end{tikzcd}$ is a $2\mhyphen crossed\mhyphen$Lie algebra, then the $g_{i,k}$'s inherit Lie algebra structures by the formulas : \begin{itemize}
    \item $[x,y]_{g_{0,0}}:=[x,y]_g$
    \item $[x,y]_{g_{0,1}}:=[x,d^v_0y]_g({}=[d^v_0x,y]_g)$
    \item $[x,y]_{g_{1,0}}:=[x,d^h_0y]_g({}=[d^h_0x,y]_g)$
    \item $[x,y]_{g_{1,1}}:=[x,d^h_0d^v_1y]_g({}=[d^h_1x,d^v_1y]_g=[d^v_1x,d^h_1y]_g= [d^v_0d^h_1x,y])$
\end{itemize} And the $2\mhyphen$fold category internal to Lie algebras has the following shape : \[\begin{tikzcd}
{(g_{1,1}\rtimes g_{0,1})\rtimes(g_{1,0}\rtimes g_{0,0})} \arrow[d, shift left=26] \arrow[r, shift left=2] \arrow[r, shift right=2] & {g_{1,0}\rtimes g_{0,0}} \arrow[d, shift right=11] \arrow[l] \\
{g_{0,1}\rtimes g_{0,0}} \arrow[r, shift left=2] \arrow[r, shift right=2] \arrow[u, shift right=31]                                                       & {g_{0,0}} \arrow[l]                                                   
\end{tikzcd}\]
\end{example}
\bibliographystyle{abbrv}
\bibliography{biblio}
\end{document}